\documentclass[11pt]{article}

\usepackage[top=50mm, bottom=50mm, left=50mm, right=50mm]{geometry}
\usepackage{lineno}                 % Line numbers
\usepackage{amssymb}                % Math
\usepackage{amsmath}                % Math
\usepackage{amsthm}                 % Math
\usepackage{mathrsfs}               % Math, \mathscr{}
\usepackage{epsfig}                 % eps figures
\usepackage{graphicx}               % Figures
\usepackage{graphics}               % ?
\usepackage{float}                  % 
\usepackage{multirow}               % 
\usepackage{color}                  %
\usepackage{fullpage}               %
\usepackage[normalem]{ulem}         % 
\usepackage{makeidx}                %
\usepackage{xspace}                 %
\usepackage{wrapfig}                %

\usepackage{url}                    % simple URL typesetting
\usepackage{booktabs}               % professional-quality tables
\usepackage{amsfonts}               % blackboard math symbols
\usepackage{nicefrac}               % compact symbols for 1/2, etc.
\usepackage{microtype}              % microtypography

\usepackage{caption}                %
\usepackage{enumitem}               %
\usepackage[export]{adjustbox}      %
\usepackage{comment}                %
\usepackage{xstring}                %
\usepackage{amsfonts}               %
\usepackage{subcaption}             %
\usepackage{tikz}                   % 
\usepackage{pgfplots}               %
\usepackage{pgfplotstable}          %

\usepackage[affil-it]{authblk}                % author and affiliations 

\usepackage[citestyle=numeric-comp,sorting=none,backend=biber]{biblatex}
\usepackage[title]{appendix}

\makeindex

\providecommand{\keywords}[1]
{
  \small	
  \textbf{\textit{Keywords: }} #1
}

\pgfplotsset{compat=1.16}

\newtheorem{problem}{Problem}
\newtheorem{dproblem}{Discrete Problem}

\newcommand\mydef{:=}

\def\bsigma{\hbox{\boldmath$\sigma$}}
\def\bepsilon{\hbox{\boldmath$\epsilon$}}
\def\bnabla{\hbox{\boldmath$\nabla$}}

\def\disp{\hbox{$u$}}
\def\bdisp{\hbox{$\mathbf{u}$}}

\def\testbdisp{\hbox{$\delta\mathbf{u}$}}

\def\strain{\hbox{$\bepsilon$}}

\def\stress{\hbox{$\bsigma$}}
\def\stresspos{\hbox{$\bsigma^+$}}
\def\stressneg{\hbox{$\bsigma^-$}}
\newcommand{\lame}[1]{%
    \IfEqCase{#1}{%
        {1}{\lambda}%
        {2}{\mu}%
    }[\PackageError{lame}{Undefined option to lame: #1}{}]%
}% Lame constants

\def\pf{\hbox{$\varphi$}}

\def\l{\hbox{$l$}}
\def\gc{\hbox{$G_c$}}

\def\mpf{\hbox{$d$}}
\def\testmpf{\hbox{$\delta\mpf$}}

\def\dom{\hbox{$\Omega$}}
\def\surf{\hbox{$\Gamma$}}

\newcommand{\surfarg}[2]{\hbox{$\surf_{#1}^{\scriptsize{#2}}$}}
\newcommand{\trac}[1]{\hbox{$\mathbf{t}_{p}^{\scriptsize{#1}}$}}

\def\surfcrack{\hbox{$\mathcal{C}$}}

\usetikzlibrary{external}
\usetikzlibrary{fadings}
\tikzfading[name=fade out, 
    inner color=transparent!0,
    outer color=transparent!100]
\usetikzlibrary{shapes}
\usetikzlibrary{decorations.pathreplacing}
\usetikzlibrary{arrows}
\usepgfplotslibrary{fillbetween}

\DeclareCaptionFont{mysize}{\fontsize{9}{9.6}\selectfont}
\begin{document}

\title{A micromorphic phase-field model for brittle and quasi-brittle fracture }

\author[1]{Ritukesh Bharali}
\author[1]{Fredrik Larsson}
\author[2]{Ralf J\"anicke}
\affil[1]{Department of Industrial and Material Science, Chalmers University of Technology}
\affil[2]{Institute of Applied Mechanics, Technische Universit\"at Braunschweig}

\date{}

\maketitle

\section*{Highlights}
\begin{itemize}
    \item \textit{First} micromorphic approach towards phase-field fracture modelling.
    \item Phase-field is transformed to a local quantity and micromorphic variable regularises the problem.
    \item Model admits point-wise (local) treatment of fracture irreversibility.
    \item System precision fracture irreversibility achieved.
    \item Dimension of the problem remains same as conventional phase-field models.
    \item Numerical experiments carried out on benchmark brittle and quasi-brittle fracture problems.
\end{itemize}

\bigskip

\begin{abstract}
\noindent The phase-field model for fracture, despite its popularity and ease of implementation comes with its set of computational challenges. They are the non-convex energy functional, variational inequality due to fracture irreversibility, the need for extremely fine meshes to resolve the fracture. In this manuscript, the focus is on the numerical treatment of variational inequality. In this context, the popular history-variable approach suffers from variationally inconsistency and non-quantifiable nature of the error introduced. A better alternative, the penalisation approach, has the potential to render the stiffness matrix ill-conditioned. In order to circumvent both aforementioned issues, a micromorphic approach towards phase-field fracture modelling is proposed in this manuscript. Within this approach, a micromorphic extension of the energy functional is carried out. This transforms the phase-field into a local variable, while introducing a micromorphic variable that regularises the fracture problem. This reduction the regularity requirements for the phase-field enables an easier implementation of the fracture irreversibility constraint through simple `max' operation, with system level precision. Numerical experiments carried out on benchmark brittle and quasi-brittle problems demonstrate the applicability and efficacy of the proposed model for a wide range of fracture problems.

% Penalising the difference between the phase-field and the micromorphic variable allows the quantification of the error introduced, unlike the history-variable approach. Through numerical experiments on benchmark problems, it is demonstrated that the micromorphic variable allows a penalty term several orders magnitude fewer than that required by a conventional penalisation technique. This results in a better conditioned stiffness matrix.    
\end{abstract}

% keywords can be removed
\keywords{phase-field fracture, brittle, quasi-brittle, micromorphic, monolithic, fracture irreversibility}

\newpage

\section{Introduction}\label{sec1}

The phase-field (smeared) fracture model is a promising alternative to conventional discrete fracture modelling techniques like XFEM \cite{sukumar2000extended,moes2002extended}, Cohesive Zone models \cite{DUGDALE1960100,BARENBLATT196255,elices2002cohesive}. Unlike the discrete fracture modelling techniques, the phase-field fracture model offers a straightforward handling of topologically complex (branching, kinking and merging of cracks) fractures based on energy minimization. Moreover, it also has the ability to operate on fixed meshes, thus circumventing the need for tedious re-meshing techniques. A comprehensive comparison between discrete and smeared fracture modelling techniques has been carried out in \cite{Egger2019}.

The phase-field fracture model emerged from the variational treatment of the Griffith fracture criterion in \cite{Francfort1998}, and its numerical adaptation in \cite{Bourdin2000,Bourdin2007}. The model introduces an auxiliary scalar variable, the phase-field, which interpolates between intact and fully broken (fractured) material states. The phase-field fracture model was cast into a thermodynamically consistent fracture in \cite{miehe2010b}. Since then, the model has been extended towards ductile fracture \cite{miehe2015486,ambati2015ductile}, anisotropic fracture \cite{TEICHTMEISTER20171,BLEYER2018213}, hydraulic fracture \cite{WILSON2016264,HEIDER201738}, desiccation cracking \cite{cajuhi2018phase,HU2020113106}, corrosion \cite{MARTINEZPANEDA2018742,KRISTENSEN2020104093}, fracture in thin films \cite{Mesgarnejad2013}, to cite a few applications. It is imperative to mention that in these applications, the phase-field component was confined to brittle fracture. In \cite{wu2017}, a unified phase-field fracture model was developed, encompassing both brittle and quasi-brittle fractures. Since then, the model has been applied in the investigation of size effect of concrete \cite{FENG201866}, hydrogen assisted cracking \cite{WU2020112614}, electro-mechanical fracture in piezo-electric solids \cite{WU2021114125} , and fracture of thermo-elastic solids \cite{MANDAL2021113648}, to cite a few. Although, comparatively fewer compared to single-scale applications, the phase-field fracture model has been studied in the context of multi-scale techniques as well. The phase-field fracture model has been extended towards concurrent multi-scale modelling in \cite{patil2019multiscale,gerasimov2018non,nguyen2019multiscale,triantafyllou2020generalized} and hierarchical multi-scale modelling in \cite{he2020numerical,Bharali2021}.

The popularity and ease of implementation of the phase-field fracture model, however, comes at the cost of its own set of computational challenges. These include (A) minimising a non-convex energy functional, (B) treatment of variational inequality due to the fracture irreversibility constraint, and (C) the need for extremely fine meshes to resolve the fracture zone. For a comprehensive review of these issues, the reader is referred to \cite{de2020numerical} and references therein. In this manuscript, the focus is solely on the numerical treatment of the fracture irreversibility constraint. In this context, \cite{Gerasimov2016,GERASIMOV2019990} opted for a simple penalisation technique, \cite{Heister2015} proposed a primal-dual active set method, while \cite{Wick2017a,wick2017modified} adopted an Augmented Lagrangian formulation based on the Moreau-Yoshida indicator function. In an alternative approach, \cite{Miehe2010} introduced an implicit history-variable as the fracture driving energy in order to ensure fracture irreversibility. Among all techniques, the penalisation approach and the history-variable approach remain popular, owing to their ease in implementation into existing finite element frameworks. However, the penalisation approach has the potential to render a stiffness matrix ill-conditioned, particularly when a stricter irreversibility tolerance is desired. This is clear from the expressions for the penalty term derived in Section 3.3.3 in \cite{GERASIMOV2019990}. The history-variable circumvents the ill-conditioning of the stiffness matrix, however, it results in the loss of variational consistency of the problem. Furthermore, to the best of the authors' knowledge, the error introduced with the history-variable approach is unknown. In order to circumvent potential ill-conditioning of stiffness matrix and non-quantifiable error, a micromorphic approach to phase-field fracture model is proposed in this manuscript. 

The micromorphic approach for gradient elasticity, viscoplasticity and damage was proposed in \cite{forest2009micromorphic}. Thereafter, the work has been extended towards crystal plasticity \cite{Forest2014,ASLAN20111311,LINDROOS2022103187}, small and finite deformation plasticity coupled with damage \cite{GRAMMENOUDIS2010140,GRAMMENOUDIS2009957} to cite a few. The micromorphic approach has been adopted in the context of ductile phase-field fracture in \cite{Miehe2016micromorphic} for introducing a length-scale to the plastic zone. To the best of the authors' knowledge, the micromorphic approach has not been explored in regard to brittle and quasi-brittle phase-field fracture models. This research gap is addressed in this manuscript, and subsequently, a simple way of enforcing the phase-field fracture irreversbility is demonstrated. To this end, the phase-field fracture energy functional is extended in the spirit of \cite{forest2009micromorphic}. This transforms the phase-field into a local quantity, while introducing a `\textit{new}' micromorphic variable that regularises the problem. The local nature of the phase-field enables a simpler treatment\footnote{in comparison to non-local quantities in the $H^1$ space} of the fracture irreversibility constraint, at material points (integration points). Furthermore, the dimension of the problem remains unchanged, with displacement and micromorphic variable being the nodal field quantities and the phase-field being a local material point variable.

The efficacy of the novel micromorphic phase-field fracture model is demonstrated with benchmark problems in both brittle and quasi-brittle fracture. These include a single edge notched specimen under tension, the Winkler L-panel experiment \cite{winkler2001traglastuntersuchungen}, and a three-point bending experiment carried out in \cite{Rots1988}. This manuscript is structured as follows: Section \ref{sec2} introduces the reader to the phase-field model for fracture, its underlying energy functional and pertinent Euler-Lagrange equations. Subsequently, in Section \ref{sec3}, the micromorphic approach for phase-field fracture is introduced. The numerical benchmark problems are addressed in Section \ref{sec4}, followed by concluding remarks in Section \ref{sec5}.

\section{Phase-field fracture model}\label{sec2}

\subsection{The energy functional}

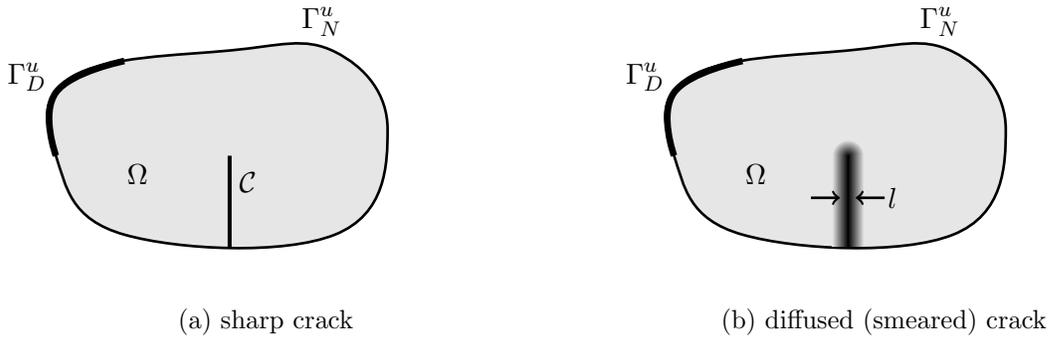
\begin{figure*}[ht!]
  \begin{subfigure}[b]{0.49\textwidth}
    \begin{tikzpicture}[scale=0.7]
    \coordinate (K) at (0,0);
    % Potato
    \draw [fill=black!10,line width=1pt] (K) plot [smooth cycle,tension=0.7] %coordinates {(3,1) (5,1.2) (7,1) (8,3) (7,4.5) (5,4.5) (2,4) (1.7,2.5)};
    coordinates {(3,1) (7,1) (8,3) (7,4.475) (5,4.5) (2,4) (1.7,2.5)};
    \node[ ] at (3.25,2.15) {$\dom$};
    % Crack Surface
    \draw[line width=1.5pt,black] (5,0.75) to (5,2.5);
    \node[ ] at (5.35,2.0) {$\surfcrack$};
    % Dirichlet Boundary
    \draw (K) [line width=2.5pt,black] plot [smooth, tension=0.8] coordinates {(3,4.3) (1.7,3.7) (1.7,2.5)};
    %\draw (K) [solid,line width=0.5pt] (1.1,4) - - (1.55,3.5);
    \node[ ] at (1.15,4) {$\surfarg{D}{u}$};
    % Neumann Boundary
    %\draw (K) [line width=2.5pt,black] plot [smooth, tension=0.7] coordinates {
    %(7.95,3.4) (7.05,4.5) (5.05,4.5)};
    %\draw (K) [solid,line width=0.5pt] (6.35,4.5) - - (5.75,4.95);
    \node[ ] at (6.75,5.05) {$\surfarg{N}{u}$};
    \end{tikzpicture}
    \caption{sharp crack}
    \label{fig:sec2:continuumpotato_a}
  \end{subfigure}
  \begin{subfigure}[b]{0.49\textwidth}
    \begin{tikzpicture}[scale=0.7]
    \coordinate (K) at (0,0);
    % Potato
    \draw [fill=black!10,line width=1pt] (K) plot [smooth cycle,tension=0.7] %coordinates {(3,1) (5,1.2) (7,1) (8,3) (7,4.5) (5,4.5) (2,4) (1.7,2.5)};
    coordinates {(3,1) (7,1) (8,3) (7,4.475) (5,4.5) (2,4) (1.7,2.5)};
    \node[ ] at (3.25,2.15) {$\dom$};
    % Crack Surface
    \fill[black, path fading=fade out, draw=none] (5,2.5) circle (0.3);
    \draw[line width=0.1pt,black] (5,0.75) to (5,2.5);
    \draw[line width=0.1pt,black] (5,0.75) to (5,2.5);
    \shade [top color=black,bottom color=black!10,shading angle=90] (5,0.78) rectangle (5.3,2.5);
    \shade [top color=black!10,bottom color=black,shading angle=90] (4.7,0.78) rectangle (5,2.5);
    \draw[->,line width=1pt,black] (4.3,1.7) to (4.85,1.7);
    \draw[<-,line width=1pt,black] (5.15,1.7) to (5.7,1.7);
    \node[ ] at (5.85,1.7) {$\l$};
    %\draw[line width=0.5pt,black] (5,2.5) to (5.5,2.8);
    %\node[ ] at (6,2.9) {$\surfpfreg$};
    % Dirichlet Boundary
    \draw (K) [line width=2.5pt,black] plot [smooth, tension=0.8] coordinates {(3,4.3) (1.7,3.7) (1.7,2.5)};
    \node[ ] at (1.15,4) {$\surfarg{D}{u}$};
    % Neumann Boundary
    %\draw (K) [line width=2.5pt,black] plot [smooth, tension=0.7] coordinates {
    %(7.95,3.4) (7.05,4.5) (5.05,4.5)};
    \node[ ] at (6.75,5.05) {$\surfarg{N}{u}$};
    \end{tikzpicture}
    \caption{diffused (smeared) crack}
    \label{fig:sec2:continuumpotato_b}
  \end{subfigure}
  \caption{A solid $\dom \in \mathbb{R}^2$ embedded with (a). sharp crack $\surfcrack$ and (b). diffused (smeared) crack, with Dirichlet and Neumann boundaries indicated as $\surfarg{D}{u}$ and $\surfarg{N}{u}$ respectively. Figure adopted from \cite{Bharali2021}}.
    \label{fig:sec2:continuumpotato}
\end{figure*}

Figure \ref{fig:sec2:continuumpotato_a} and \ref{fig:sec2:continuumpotato_b} illustrates a discrete fracture and a phase-field regularised fracture in a 2D continuum respectively. The salient feature of the a discrete fracture model is the explicit representation of the fracture through displacement discontinuities. The phase-field fracture model belongs to the smeared representation of fracture, where the phase-field, $\pf \in [0,1]$ interpolates between intact and fully broken (fractured) material states. The fractured continuum in both figures is assumed to occupy a domain $\dom \in \mathbb{R}^{\text{dim}}$ ($\text{dim} = 2$ in this case). The external boundary $\surf$ comprises of the Dirichlet and Neumann boundaries, represented by $\surfarg{D}{u}$ and $\surfarg{N}{u}$ respectively. Note that $\surf = \surfarg{D}{u} \cup \surfarg{N}{u}$ and $\surfarg{D}{u} \cap \surfarg{N}{u} = \emptyset$. 

The energy functional pertaining to the phase-field fracture model \cite{de2020numerical} is given by,

\begin{equation}{\label{eqn:sec2:EFunc}}
    \displaystyle E(\bdisp,\pf) = \int_{\dom}^{} g(\pf)  \Psi^+(\strain[\bdisp]) \: \text{d}\dom + \int_{\dom}^{} \Psi^-(\strain[\bdisp]) \: \text{d}\dom - \int_{\surfarg{N}{u}}^{}  \trac{u} \cdot \bdisp \: \text{d}\surf + \int_{\dom}^{} \dfrac{\gc}{c_w \l} \left( w(\pf) + \l^2 |\bnabla \pf|^2 \right) \: \text{d}\dom,
\end{equation}

where, a degradation function $g(\pf)$ is attached to the fracture driving strain energy density $\Psi^+$, with $\Psi^-$ and $\trac{u}$ being the residual strain energy density and the prescribed traction respectively. The last integral in the above equation represents the fracture energy, where $\gc$ and $\l$ are the Griffith fracture energy and the fracture length-scale respectively, while $c_w$ is a normalisation constant associated with the choice of the locally dissipated fracture energy function $w(\pf)$. Some commonly adopted phase-field fracture models are presented in the Table \ref{sec2:table:brittleQuasiModelsParams}. 

\begingroup
\renewcommand{\arraystretch}{1.175}
\begin{table}[ht!]
    \centering
    \begin{tabular}{lllc} \hline
    Fracture model  & $w(\pf)$ & $c_w$ & $g(\pf)$  \\ \hline
    Brittle AT1 \, \cite{pham2011} & $\pf$  & 8/3 & $(1-\pf)^2$ \\ 
    Brittle AT2 \, \cite{Bourdin2000}  & $\pf^2$ & 2 &  $(1-\pf)^2$ \\ 
    Quasi-brittle \cite{wu2017}  & $2\pf-\pf^2$ & $\pi$ & $\dfrac{(1-\pf)^p}{(1-\pf)^p + a_1\pf + a_1 a_2 \pf^2 + a_1 a_2 a_3 \pf^3}$ \\ \hline
    \end{tabular}
    \caption{Brittle and quasi-brittle phase-field fracture models.}
    \label{sec2:table:brittleQuasiModelsParams}
\end{table}
\endgroup

The quasi-brittle phase-field fracture model requires additional parameters $p$, $a_1$, $a_2$ and $a_3$, as observed from the Table \ref{sec2:table:brittleQuasiModelsParams}. These parameters are chosen in way, such that different traction-separation laws may be obtained. For a detailed derivation of these parameters, the reader is referred to \cite{wu2017}. Following \cite{wu2017}, the constant $a_1$ is given by,

\begin{equation}
    a_1 = \dfrac{4 E_0 \gc}{\pi \l f_t^2},
\end{equation}

\noindent where, the newly introduced parameters $E_0$ and $f_t$ represent the Young's Modulus and the material tensile strength respectively. The other parameters vary for different traction-separation laws, and given in Table \ref{sec2:table:quasiBrittleParams}.

\begingroup
\renewcommand{\arraystretch}{1.175}
\begin{table}[ht!]
    \centering
    \begin{tabular}{llll} \hline
     Traction-separation law & $p$  & $a_2$  & $a_3$ \\ \hline
    Linear Softening   & $2$  & $-0.5$  & $0$ \\ 
    Exponential Softening   & $2.5$  & $2^{5/3}-3$  & $0$ \\ 
    Cornelissen et. al.\cite{cornelissen1986experimental} Softening  & $2$  & $1.3868$  & $0.6567$ \\ \hline
    \end{tabular}
    \caption{Quasi-brittle phase-field fracture model parameters for different traction-separation laws \cite{wu2017}.}
    \label{sec2:table:quasiBrittleParams}
\end{table}
\endgroup

As demonstrated in Equation (\ref{eqn:sec2:EFunc}), the strain energy density is split into fracture driving and residual components, $\Psi^+$ and $\Psi^-$. In this manuscript, the common consensus that a fracture occurs under tensile loading is adopted. Therefore, the fracture driving and residual strain energy densities are computed using the Amor split \cite{pham2011}. With this approach, the fracture is driven by the deviatoric part of the strain energy density together with a positive volumetric part of the strain energy density. The negative volumetric part corresponds to the residual energy. Mathematically,

\begin{equation}
    \Psi^+ = \frac{1}{2} K \langle tr(\strain) \rangle_+^2 + \mu ( \strain_{dev} \colon \strain_{dev} )
\end{equation}

and

\begin{equation}
    \Psi^+- = \frac{1}{2} K \langle tr(\strain) \rangle_-^2.
\end{equation}

\noindent Here, $K$ and $\mu$ represents the bulk and the shear modulus of the material respectively, $tr$ is a trace operator, $\langle\bullet\rangle_{\pm}$ represents the positive/negative Macaulay brackets, and $\epsilon_{dev}$ is the deviatoric strain. The corresponding Cauchy stresses are given by,

\begin{equation}
    \stress^+ \mydef \dfrac{\partial \Psi^+}{\partial \strain} = K \mathcal{H} \big( tr(\strain) \big) \, \mathbf{I} + 2\mu\strain_{dev}
\end{equation}

and 

\begin{equation}
    \stress^- \mydef \dfrac{\partial \Psi^-}{\partial \strain} = K \mathcal{H} \big( - tr(\strain) \big) \, \mathbf{I},
\end{equation}

\noindent adopting the heaviside function $\mathcal{H}$.

\subsection{Euler-Lagrange equations}

In order to simulate fracture initiation and propagation in continuum $\dom$, minimisation of the energy functional (\ref{eqn:sec2:EFunc}) is required, w.r.t. its solution variables $\bdisp$ and $\pf$. This results in the Euler-Lagrange equations. Along with appropriately defined test and trial spaces, they result in the following problem: 

\begin{problem}\label{Problem1}
Find ($\bdisp$, $\pf$) $\in \mathbb{U} \times \mathbb{P}$ with

\begin{subequations}
\begin{align}
\int_{\dom}^{} g(\pf) \stresspos \colon\strain[\testbdisp] \: \normalfont \text{d} \dom + \int_{\dom}^{} \stressneg \colon \strain[\testbdisp] \: \text{d}\dom 
- \int_{\surfarg{N}{\bdisp}}^{} \trac{\disp} \: \testbdisp \: \normalfont \text{d}\surf &= 0  & \: \forall \: \testbdisp \in \mathbb{U}^0, \label{eqn:sec2:EL_momentumbalance} \\
\int_{\dom}^{} \dfrac{\gc}{c_w \l} \left( w'(\pf) \, (\hat{\pf}-\pf) + 2\l^2 \, \bnabla \pf\cdot \bnabla (\hat{\pf}-\pf) \right) \normalfont \text{d}\dom + \int_{\dom}^{} g'(\pf) \Psi^+ (\hat{\pf}-\pf) \: \text{d}\dom & \geq 0 & \: \forall \: \hat{\pf} \in \mathbb{P}, \label{eqn:sec2:pf_evolution} 
\end{align}
\end{subequations}
\noindent using pertinent time-dependent Dirichlet boundary conditions $\bdisp^p$ on $\surfarg{D}{u}$ and $\pf^p$ on $\surfarg{D}{\pf}$, and and Neumann boundary condition $\trac{\bdisp}$ on $\surfarg{N}{\bdisp}$. The trial and test spaces are defined as
\begin{subequations}
\begin{align}
\mathbb{U} & = \{ \bdisp \in [H^1(\dom)]^{\normalfont\text{dim}}| \bdisp = \bdisp^p \text{ on } \surfarg{D}{u} \},  {\label{eqn:sec2:disp_space}} \\ 
\mathbb{U}^0 & = \{ \bdisp \in [H^1(\dom)]^{\normalfont\text{dim}}| \bdisp = \mathbf{0} \text{ on } \surfarg{D}{u} \},  {\label{eqn:sec2:test_disp_space}} \\
\mathbb{P} & = \{ \pf \in [H^1(\dom)]^1 | \pf \geq {}^{n}\pf |\pf = \pf^p \text{ on } \surfarg{D}{\pf} \}.  {\label{eqn:sec2:pf_space}}
\end{align}
\end{subequations}
Note that the requirement $\pf \geq {}^{n}\pf$ in (\ref{eqn:sec2:pf_space}) ensures fracture irreversibility, with $n$ referring to the previous time-step. {\color{black}\hfill $\blacksquare$}
\end{problem}

From Problem \ref{Problem1}, it is clear that the fracture irreversibility constraint manifests in the form of a variationally inequality Euler-Lagrange equation (\ref{eqn:sec2:pf_evolution}) with restrictive trial and test set (\ref{eqn:sec2:pf_space}) for the phase-field. As mentioned in Section \ref{sec1}, several researchers have proposed different methods to treat the variational inequality problem. In the next section, a micromorphic phase-field fracture model is proposed, where the phase-field variable becomes a local quantity, while a micromorphic variable ensures regularisation of the problem. This enables a simpler `\textit{point-wise}' treatment of the fracture irreversibility constraint.

\section{Micromorphic phase-field fracture model}\label{sec3}

\subsection{Energy functional and variants}

In this sub-section, a micromorphic extension of the phase-field fracture energy functional (\ref{eqn:sec2:EFunc}) is carried out in the spirit of \cite{forest2009micromorphic}. This results in,

\begin{align}
\Tilde{E}(\bdisp,\pf,\mpf) & = \int_{\dom}^{} g(\pf) \Psi^+(\strain[\bdisp]) \: \text{d}\dom + \int_{\dom}^{} \Psi^-(\strain[\bdisp]) \: \text{d}\dom - \int_{\surfarg{N}{u}}^{}  \trac{u} \cdot \bdisp \: \text{d}\surf   {\label{eqn:sec3:microE1}} \\ 
& + \int_{\dom}^{} \dfrac{\gc}{c_w \l} \left( w(\pf) + \l^2 |\bnabla \mpf|^2 \right) \: \text{d}\dom + \int_{\dom}^{} \frac{\alpha}{2} (\pf-\mpf)^2 \: \text{d}\dom. \nonumber
\end{align}

Here, $\mpf$ is the `\textit{new}' micromorphic variable that regularises the fracture problem. From the above equation, it is observed that the micromorphic variable $\mpf$ is solely associated with the gradient term in the fracture energy integral, and the phase-field $\pf$ becomes a local quantity. Consequently, the regularity requirements on the phase-field w.r.t. the existence of derivatives is circumvented. Moreover, the micromorphic approach also introduces an \textit{additional energy} term associated with the difference between the phase-field and the micromorphic variable. Theoretically, for the continuous problem, in the limit, the interaction parameter, $\alpha \rightarrow \infty$, the original energy functional (\ref{eqn:sec2:EFunc}) is recovered.

\subsection{Euler-Lagrange equations}

The set of Euler-Lagrange equations for the micromorphic phase-field fracture model is obtained upon minimising the energy functional (\ref{eqn:sec3:microE1}) w.r.t. its solution variables $\bdisp$, $\pf$ and $\mpf$. Along with appropriately defined test and trial spaces, it results in the following problem:

\begin{problem}\label{Problem2}
Find ($\bdisp$, $\pf$, $\mpf$) $\in \mathbb{U} \times \mathbb{P} \times \mathbb{D}$ with

\begin{subequations}
\begin{align}
\int_{\dom}^{} g(\pf) \stresspos \colon\strain[\testbdisp] \: \normalfont \text{d} \dom + \int_{\dom}^{} \stressneg \colon \strain[\testbdisp] \: \text{d}\dom - \int_{\surfarg{N}{\bdisp}}^{} \trac{\disp} \: \testbdisp \: \normalfont \text{d}\surf & = 0  & \: \forall \: \testbdisp \in \mathbb{U}^0, \label{eqn:P2:EL_momentumbalance} \\
\int_{\dom}^{} g'(\pf) \Psi^+ (\hat{\pf}-\pf) \normalfont \: \text{d}\dom + \int_{\dom}^{} \dfrac{\gc}{c_w \l} w'(\pf) (\hat{\pf}-\pf) \: \text{d}\dom + \int_{\dom}^{} \alpha (\pf-\mpf) (\hat{\pf}-\pf) \: \text{d}\dom & \geq 0 & \: \forall \: \hat{\pf} \in \mathbb{P}, \label{eqn:P2:pf_evolution} \\
\int_{\dom}^{} \dfrac{2 \gc \l}{c_w} \, \bnabla \mpf\cdot \bnabla \testmpf \normalfont \: \text{d}\dom - \int_{\dom}^{} \alpha (\pf-\mpf) \testmpf \: \text{d}\dom & = 0 & \: \forall \: \testmpf \in \mathbb{D}, \label{eqn:P2:mpf_evolution}
\end{align}
\end{subequations}

\noindent using pertinent time-dependent Dirichlet and Neumann boundary conditions, $\bdisp^p$ on $\surfarg{D}{u}$, and $\trac{\bdisp}$ on $\surfarg{N}{\bdisp}$ respectively. The trial and test spaces are given by,

\begin{subequations}
\begin{align}
\mathbb{U} & = \{ \bdisp \in [H^1(\dom)]^{\normalfont\text{dim}}| \bdisp = \bdisp^p \text{ on } \surfarg{D}{u} \},  {\label{eqn:P2:disp_space}} \\ 
\mathbb{U}^0 & = \{ \bdisp \in [H^1(\dom)]^{\normalfont\text{dim}}| \bdisp = \mathbf{0} \text{ on } \surfarg{D}{u} \},  {\label{eqn:P2:test_disp_space}} \\
\mathbb{D} & = \{ \mpf \in [H^1(\dom)]^1] \},  {\label{eqn:P2:mpf_space}} \\ 
\mathbb{P} & = \{ \pf \in [L^2(\dom)] | \pf \geq {}^{n}\pf \}.  {\label{eqn:P2:pf_space}}
\end{align}
\end{subequations}

{\color{black}\hfill $\blacksquare$}
\end{problem}

In Problem \ref{Problem2}, the phase-field evolution equation (\ref{eqn:P2:pf_evolution}) is local. As such, $\pf$ can be computed `point-wise' in the computational domain. In particular, the root(s) of the possibly nonlinear scalar equation,

\begin{equation}\label{eqn:sec:2-4:localPfEqn}
    g'(\pf) \Psi^+(\strain[\bdisp]) + \dfrac{\gc}{c_w \l} w'(\pf) + \alpha (\pf - \mpf) = 0
\end{equation}

\noindent yield(s) the phase-field $\pf$, if ${}^{n}\pf < \pf < 1$. Using the locally computed phase-field, the complete problem for micromorphic phase-field fracture is assumes the form:

\begin{problem}\label{Problem3}
Find ($\bdisp$, $\mpf$) $\in \mathbb{U} \times \mathbb{D}$ with

\begin{subequations}
\begin{align}
\int_{\dom}^{} g(\pf) \stresspos \colon\strain[\testbdisp] \: \normalfont \text{d} \dom + \int_{\dom}^{} \stressneg \colon \strain[\testbdisp] \: \text{d}\dom - \int_{\surfarg{N}{\bdisp}}^{} \trac{\disp} \: \testbdisp \: \normalfont \text{d}\surf & = 0  & \: \forall \: \testbdisp \in \mathbb{U}^0, \label{eqn:P3:EL_momentumbalance} \\
\int_{\dom}^{} \dfrac{2 \gc \l}{c_w} \, \bnabla \mpf\cdot \bnabla \testmpf \normalfont \: \text{d}\dom - \int_{\dom}^{} \alpha (\pf-\mpf) \testmpf \: \text{d}\dom & = 0 & \: \forall \: \testmpf \in \mathbb{D}, \label{eqn:P3:mpf_evolution}
\end{align}
\end{subequations}

\noindent using pertinent time-dependent Dirichlet boundary conditions $\bdisp^p$ on $\surfarg{D}{u}$, and Neumann boundary condition $\trac{u}$ on $\surfarg{N}{u}$. The trial and test spaces are defined as
\begin{subequations}
\begin{align}
\mathbb{U} & = \{ \bdisp \in [H^1(\dom)]^{\normalfont\text{dim}}| \bdisp = \bdisp^p \text{ on } \surfarg{D}{u} \},  {\label{eqn:P3:disp_trialspace}} \\ 
\mathbb{U}^0 & = \{ \bdisp \in [H^1(\dom)]^{\normalfont\text{dim}}| \bdisp = \mathbf{0} \text{ on } \surfarg{D}{u} \},  {\label{eqn:P3:disp_testspace}} \\
\mathbb{D} & = \{ \mpf \in [H^1(\dom)] \},  {\label{eqn:P3:trialspace}}
\end{align}
\end{subequations} 

with the local phase-field $\pf$ computed using (\ref{eqn:sec:2-4:localPfEqn}). {\color{black}\hfill $\blacksquare$}
\end{problem}

It is worth mentioning that the local phase-field evolution (\ref{eqn:sec:2-4:localPfEqn}) is linear for brittle AT1 and AT2 fracture models. This is due to the quadratic nature of the degradation function $g(\pf) = (1-\pf)^2$ coupled with linear/quadratic locally dissipated fracture energy function $w(\pf) = \pf$ and $\pf^2$ for the AT1 and AT2 models respectively. This yield explicit expressions for the local phase-field variable,

\begin{equation}\label{eqn:sec3-2:explicitPfAT1}
    \pf = \operatorname{min} \bigg( \operatorname{max} \Bigg( \dfrac{2 \Psi^+ + \alpha \mpf - \frac{3 \gc}{8 \l}}{2 \Psi^+ + \alpha}, {}^{n}\pf \Bigg), 1 \bigg),
\end{equation}

\noindent for AT1, and

\begin{equation}\label{eqn:sec3-2:explicitPfAT2}
    \pf = \operatorname{min} \bigg( \operatorname{max} \Bigg( \dfrac{2 \Psi^+ + \alpha \mpf}{2 \Psi^+ + \alpha + \frac{\gc}{\l}}, {}^{n}\pf \Bigg), 1 \bigg).
\end{equation}

\noindent for AT2 model respectively. In the case of quasi-brittle phase-field fracture model, the rational nature of the degradation function (see Table \ref{sec2:table:brittleQuasiModelsParams}) results in nonlinear scalar equation for the local phase-field variable, warranting the need for the Newton-Raphson method.

\subsection{Discrete equations}

The work in this manuscript is executed within the framework of the finite element method \cite{zienkiewicz1977finite,hughes2012finite}, using triangular (T3) elements for both, displacement and micromorphic fields. This allows assuming the displacement and the micromorphic fields at the nodes ($\tilde{\bdisp}_i,\tilde{\mpf}_i$) as the primary unknowns, with the corresponding continuous fields ($\bdisp$,$\mpf$) approximated as, 

\begin{equation}\label{eqn:discrete_disp_mpf}
  \bdisp = \sum\limits_{i = 1}^m {N_i^{\boldsymbol{u}}\tilde{\boldsymbol{u}}_i}\;\;\;\text{ , }
\mpf = \sum\limits_{i = 1}^m N_i^{d}\tilde{\mpf}_i.
\end{equation}

\noindent In the above equation, $N_i^{\boldsymbol{u}}$ and $N_i^{d}$ are the interpolation functions for the displacement and the micromorphic phase-field, associated with the $i^{\text{th}}$ node. The spatial derivatives of the interpolation functions $N_i^{\boldsymbol{u}}$ and $N_i^{\mpf}$ in a two-dimensional case are given by,

\begin{equation}\label{eqn:Bmatrices}
    \mathbf{B}_i^{\boldsymbol{u}} = \left[ {\begin{array}{*{20}{c}}
  {\begin{array}{*{20}{c}}
  {{N_{i,x}}} \\ 
  0 \\ 
  {{N_{i,y}}}\\
\end{array}}&{\begin{array}{*{20}{c}}
  0 \\ 
  {{N_{i,y}}} \\ 
  {{N_{i,x}}} \\
\end{array}}
\end{array}} \right] \text{ , } 
\mathbf{B}_i^{d}  = \left[ {\begin{array}{*{20}{c}}
  {{N_{i,x}}} \\ 
  {{N_{i,y}}}
\end{array}} \right].
\end{equation}

\noindent Here, the subscripts $,x$ and $,y$ indicate spatial derivatives in $x$ and $y$ directions respectively. Using (\ref{eqn:Bmatrices}), the strain $\strain$, and the gradient of the micromorphic variable $\bnabla\mpf$ are defined as,

\begin{equation}\label{eqn:strain_gradpf}
  \boldsymbol{\epsilon} = \sum\limits_{i = 1}^m {\mathbf{B}_i^{\boldsymbol{u}} \tilde{\boldsymbol{u}}_i}\;\;\;\text{ , }
\bnabla\mpf = \sum\limits_{i = 1}^m {\mathbf{B}_i^{d} \tilde{\mpf}}.
\end{equation}

\noindent The discrete phase-field fracture problem is obtained upon inserting (\ref{eqn:discrete_disp_mpf}-\ref{eqn:strain_gradpf}) in the Euler-Lagrange equations from Problem \ref{Problem3}. Thereafter, (\ref{eqn:P3:EL_momentumbalance}) and (\ref{eqn:P3:mpf_evolution}) are assumed as the internal forces, and stiffness matrix derived from its derivative. This notation is consistent with \cite{de2012nonlinear}, and allows the presentation of the phase-field fracture problem in the incremental iterative framework as:

\smallskip
\begin{dproblem}\label{Problem4}
Compute the solution increment ($\Delta\tilde{\bdisp}$, $\Delta\tilde{\mpf}$)$_{i+1}$ in the current iteration $i+1$ using

\begin{subequations}
\begin{equation}\label{dproblem1:discSystem}
\underbrace{
\begin{bmatrix}
\mathbf{K}^{\boldsymbol{uu}} & \mathbf{K}^{\boldsymbol{u}d} \\ 
\mathbf{K}^{d\boldsymbol{u}} & \mathbf{K}^{dd}
\end{bmatrix}_{i}}_{\text{Stiffness matrix}} \begin{Bmatrix}
\Delta \tilde{\bdisp} \\ 
\Delta \tilde{\mpf}
\end{Bmatrix}_{i+1}  =  \underbrace{\begin{Bmatrix}
\mathbf{f}^{ext,\boldsymbol{u}} \\ 
\mathbf{f}^{ext,d}
\end{Bmatrix}_{i} - \begin{Bmatrix}
\mathbf{f}^{int,\boldsymbol{u}} \\ 
\mathbf{f}^{int,d}
\end{Bmatrix}_{i}}_{\text{Residual}},
\end{equation}
\noindent and update the solution fields,
\begin{equation}\label{dproblem1:solUpdate}
\begin{Bmatrix}
\tilde{\bdisp} \\ 
\tilde{\mpf}
\end{Bmatrix}_{i+1} = \begin{Bmatrix}
\tilde{\bdisp} \\ 
\tilde{\mpf}
\end{Bmatrix}_{i} +
\begin{Bmatrix}
\Delta \tilde{\bdisp} \\ 
\Delta \tilde{\mpf}
\end{Bmatrix}_{i+1},    
\end{equation}
\noindent until a certain convergence measure is fulfilled. The local element stiffness matrices are computed as:

\begin{equation}\label{eq:stiff_comp}
\displaystyle
    \begin{split}
        \mathbf{K}^{\boldsymbol{uu}} & = \int_{\Omega} \left[\mathbf{B}^{\boldsymbol{u}}\right]^T \bigg( \underbrace{g(\pf) \dfrac{\partial \stress^+}{\partial \strain} + \dfrac{\partial \stress^-}{\partial \strain}}_{\mathbf{D}} + g'(\pf) \dfrac{\partial \pf}{\partial \strain} \stress^+ \bigg) \left[\mathbf{B}^{\boldsymbol{u}}\right]\,d\Omega, \\
        \mathbf{K}^{\boldsymbol{u}d} & =  \int_{\Omega} \left[\mathbf{B}^{\boldsymbol{u}}\right]^T \bigg( g'(\pf) \dfrac{\partial \pf}{\partial \mpf} \stress^+ \bigg) \left[N^{d} \right] \,d\Omega, \\
        \mathbf{K}^{d\boldsymbol{u}} & 
        = - \int_{\Omega} \left[N^{d}\right]^T \bigg( \alpha \dfrac{\partial \pf}{\partial \strain} \bigg) \left[\mathbf{B}^{\boldsymbol u}\right]\;d\Omega, \\
        \mathbf{K}^{dd} & = 
        \int_{\Omega}\left\{\left[\mathbf{B}^{d}\right]^T \left( \dfrac{2 \gc \l}{c_w} \right) \left[\mathbf{B}^{d}\right] + \left[N^{d}\right]^T \alpha \left( 1 - \dfrac{\partial \pf}{\partial \mpf} \right) \left[N^{d}\right] \right\} \;d\Omega, \\
    \end{split}
\end{equation}
\noindent and the local internal force vectors are computed as
\begin{equation}\label{eq:fint_comp}
    \begin{split}
        \mathbf{f}^{int,\boldsymbol{u}} & = \int_{\Omega} \left[\mathbf{B}^{\boldsymbol{u}}\right]^T \big( g(\pf) \stress^+ + \stress^- \big) \,d\Omega,\\
        \mathbf{f}^{int,d} & = 
        \int_{\Omega}\left\{\left[\mathbf{B}^{d}\right]^T \left( \dfrac{2 \gc \l}{c_w} \right) \left[\mathbf{B}^{d}\right] \tilde{\mpf} - \left[N^{d}\right]^T \alpha \left( \pf - \mpf \right) \right\} \;d\Omega.
    \end{split}
\end{equation}
\noindent The external force vectors $\mathbf{f}^{int,\boldsymbol{u}}$ and $\mathbf{f}^{int,d}$ are considered equal to zero. Computation of the material stiffness $\normalfont \mathbf{D}$ in Voigt notation is presented in Appendix \ref{appendixA}. {\color{black}\hfill $\blacksquare$}
\end{subequations}
\end{dproblem}

 The discrete problem \ref{Problem4} is non-convex. As such, conventional incremental solution techniques, like the Newton-Raphson method fails to achieve convergence in the softening regime, possibly due to an indefinite stiffness matrix. Efforts to circumvent this issue include the novel line search technique proposed in \cite{Gerasimov2016}, the use of arc-length solvers \cite{vignollet2014phase,may2015numerical,BHARALI2022114927}, trust-region methods \cite{kopanivcakova2020recursive} and convexification via extrapolation of the phase-field for the momentum balance equation \cite{Heister2015}. These aforementioned techniques are within the framework of monolithic solution techniques. In \cite{Miehe2010}, the alternate minimisation solution techniques was proposed, since the phase-field energy functional is convex with respect to either of the coupled field, if the other one is held constant. In this manuscript, a monolithic solution technique is desired. As such, the convexification strategy proposed by \cite{Heister2015} is adopted. With this choice, the phase-field for the momentum balance equation $\hat{\pf}$, is computed using,

\begin{equation}\label{eqn:sec:2-4:localExPfEqn}
    g'(\hat{\pf}) \Psi^+(\strain[\bdisp]) + \dfrac{\gc}{c_w \l} w'(\hat{\pf}) + \alpha (\hat{\pf} - \hat{\mpf}) = 0,
\end{equation}

\noindent with $\hat{\mpf}$ being the linearly extrapolated micromorphic variable from the last two converged (time) steps. It is important to note that this extrapolation-based computation of the phase-field applies only to the momemtum balance equation and not the micromorphic variable evolution equation. For the micromorphic variable evolution equation, the phase-field is computed using (\ref{eqn:sec:2-4:localPfEqn}), restated below for clarity,

\begin{equation}\label{eqn:sec:2-4:localPfEqnrepeat}
    g'(\pf) \Psi^+(\strain[\bdisp]) + \dfrac{\gc}{c_w \l} w'(\pf) + \alpha (\pf - \mpf) = 0.
\end{equation}

Consequently, the discrete problem assumes the form:

\smallskip
\begin{dproblem}\label{Problem5}
Compute the solution increment ($\Delta\tilde{\bdisp}$, $\Delta\tilde{\mpf}$)$_{i+1}$ in the current iteration $i+1$ using

\begin{subequations}
\begin{equation}\label{dproblem2:discExSystem}
\underbrace{
\begin{bmatrix}
\mathbf{K}^{\boldsymbol{uu}} & \mathbf{K}^{\boldsymbol{u}d} \\ 
\mathbf{K}^{d\boldsymbol{u}} & \mathbf{K}^{dd}
\end{bmatrix}_{i}}_{\text{Stiffness matrix}} \begin{Bmatrix}
\Delta \tilde{\bdisp} \\ 
\Delta \tilde{\mpf}
\end{Bmatrix}_{i+1}  =  \underbrace{\begin{Bmatrix}
\mathbf{f}^{ext,\boldsymbol{u}} \\ 
\mathbf{f}^{ext,d}
\end{Bmatrix}_{i} - \begin{Bmatrix}
\mathbf{f}^{int,\boldsymbol{u}} \\ 
\mathbf{f}^{int,d}
\end{Bmatrix}_{i}}_{\text{Residual}},
\end{equation}
\noindent and update the solution fields,
\begin{equation}\label{dproblem1:solUpdate}
\begin{Bmatrix}
\tilde{\bdisp} \\ 
\tilde{\mpf}
\end{Bmatrix}_{i+1} = \begin{Bmatrix}
\tilde{\bdisp} \\ 
\tilde{\mpf}
\end{Bmatrix}_{i} +
\begin{Bmatrix}
\Delta \tilde{\bdisp} \\ 
\Delta \tilde{\mpf}
\end{Bmatrix}_{i+1},    
\end{equation}
\noindent until a certain convergence measure is fulfilled. The local element stiffness matrices are computed as:

\begin{equation}\label{eq:stiff_comp2}
\displaystyle
    \begin{split}
        \mathbf{K}^{\boldsymbol{uu}} & = \int_{\Omega} \left[\mathbf{B}^{\boldsymbol{u}}\right]^T \bigg( \underbrace{g(\hat{\pf}) \dfrac{\partial \stress^+}{\partial \strain} + \dfrac{\partial \stress^-}{\partial \strain}}_{\mathbf{D}} \bigg) \left[\mathbf{B}^{\boldsymbol{u}}\right]\,d\Omega, \\
        \mathbf{K}^{\boldsymbol{u}d} & =  \mathbf{0}, \\
        \mathbf{K}^{d\boldsymbol{u}} & 
        = - \int_{\Omega} \left[N^{d}\right]^T \bigg( \alpha \dfrac{\partial \pf}{\partial \strain} \bigg) \left[\mathbf{B}^{\boldsymbol u}\right]\;d\Omega, \\
        \mathbf{K}^{dd} & = 
        \int_{\Omega}\left\{\left[\mathbf{B}^{d}\right]^T \left( \dfrac{2 \gc \l}{c_w} \right) \left[\mathbf{B}^{d}\right] + \left[N^{d}\right]^T \alpha \left( 1 - \dfrac{\partial \pf}{\partial \mpf} \right) \left[N^{d}\right] \right\} \;d\Omega, \\
    \end{split}
\end{equation}
\noindent and the local internal force vectors are computed as
\begin{equation}\label{eq:fint_comp}
    \begin{split}
        \mathbf{f}^{int,\boldsymbol{u}} & = \int_{\Omega} \left[\mathbf{B}^{\boldsymbol{u}}\right]^T \big( g(\hat{\pf}) \stress^+ + \stress^- \big) \,d\Omega,\\
        \mathbf{f}^{int,d} & = 
        \int_{\Omega}\left\{\left[\mathbf{B}^{d}\right]^T \left( \dfrac{2 \gc \l}{c_w} \right) \left[\mathbf{B}^{d}\right] \tilde{\mpf} - \left[N^{d}\right]^T \alpha \left( \pf - \mpf \right) \right\} \;d\Omega.
    \end{split}
\end{equation}
The local phase-fields $\hat{\pf}$ and $\pf$ are computed using (\ref{eqn:sec:2-4:localExPfEqn}) and (\ref{eqn:sec:2-4:localPfEqnrepeat}) respectively. The external force vectors $\mathbf{f}^{int,\boldsymbol{u}}$ and $\mathbf{f}^{int,d}$ are considered equal to zero. Computation of the material stiffness $\normalfont \mathbf{D}$ in Voigt notation is presented in Appendix \ref{appendixA}. {\color{black}\hfill $\blacksquare$}
\end{subequations}
\end{dproblem}

\section{Numerical Study}\label{sec4}

In this section, numerical experiments are carried out on brittle and quasi-brittle phase-field fracture problems. The problems in brittle fracture include the single edge notched specimen under tension and shear, and the three-point bending test \cite{Miehe2010}, whereas for quasi-brittle fracture, the Winkler concrete L-panel experiment \cite{winkler2001traglastuntersuchungen} and the concrete three point bending experiment \cite{Rots1988} are analysed. For each problem, the geometry, loading conditions as well as the additional model parameters are presented in the respective sub-sections. The phase-field fracture topology in the final step of the analysis, and the load-displacement curves are also shown therein.

All problems are solved in a fully coupled (monolithic) sense, adopting the Newton-Raphson method. The iterative procedure is terminated when an error measure defined as ratio of the norm of the residual in the current iteration to that of the first iteration is less than a certain tolerance, $tol$. In this manuscript, $tol = 1e-3$. The linear problem within each iteration is solved using the GMRES solver, preconditioned with an incomplete LU preconditioner. Moreover, for all numerical experiments, the interaction parameter $\alpha$ is parametrised as

\begin{equation}
\alpha = \beta \dfrac{\gc}{\l},
\end{equation}

with $\beta$ being a user-defined non-dimensional scalar.

\subsection{Single Edge Notched specimen under Tension (SENT)}\label{sec4:SENtension}

The single edge notched specimen \cite{Miehe2010} has been studied extensively under tensile and shear loading in the phase-field fracture literature. The geometry consists of a unit square (in mm) embedded with a horizontal notch, midway along height and equal to half of the edge length as shown in Figure \ref{sec4:fig:tension}. The notch is modelled explicitly in the finite element mesh. A quasi-static loading is applied at the top boundary in the form of prescribed displacement increment $\Delta\disp = 1e-4$[mm] for the first $55$ steps, following which it is changed to $1e-6$[mm]. The bottom boundary remains fixed. The additional model parameters are presented in Table \ref{sec4:table:miehe_tension_shear}.

\begin{figure}[ht]
\begin{minipage}[b]{0.33\linewidth}
\centering
\begin{tikzpicture}[scale=0.7]
    \coordinate (K) at (0,0);
    % Square
    \draw[line width=0.75pt,black] (-2.5,-2.5) to (2.5,-2.5);
    \draw[line width=0.75pt,black] (2.5,-2.5) to (2.5,2.5);
    \draw[line width=0.75pt,black] (2.5,2.5) to (-2.5,2.5);
    \draw[line width=0.75pt,black] (-2.5,2.5) to (-2.5,-2.5);
    \draw[line width=1.5pt,red] (-2.5,0) to (0,0);
    % Top Boundary
    \draw[line width=0.75pt,black] (2.5,2.75) to (-2.5,2.75);
    \draw[->,line width=1.5pt,black] (0.0,3.0) to (0.0,3.5);
    \node[ ] at (0,3.85) {$\Delta\disp$};
    % Bottom Boundary
    \draw[line width=0.75pt,black] (-2.5,-2.5) to (-2.75,-2.75);
    \draw[line width=0.75pt,black] (-2.25,-2.5) to (-2.5,-2.75);
    \draw[line width=0.75pt,black] (-2.0,-2.5) to (-2.25,-2.75);
    \draw[line width=0.75pt,black] (-1.75,-2.5) to (-2.0,-2.75);
    \draw[line width=0.75pt,black] (-1.5,-2.5) to (-1.75,-2.75);
    \draw[line width=0.75pt,black] (-1.25,-2.5) to (-1.5,-2.75);
    \draw[line width=0.75pt,black] (-1.0,-2.5) to (-1.25,-2.75);
    \draw[line width=0.75pt,black] (-0.75,-2.5) to (-1.0,-2.75);
    \draw[line width=0.75pt,black] (-0.5,-2.5) to (-0.75,-2.75);
    \draw[line width=0.75pt,black] (-0.25,-2.5) to (-0.5,-2.75);
    \draw[line width=0.75pt,black] (0.0,-2.5) to (-0.25,-2.75);
    \draw[line width=0.75pt,black] (0.25,-2.5) to (0.0,-2.75);
    \draw[line width=0.75pt,black] (0.5,-2.5) to (0.25,-2.75);
    \draw[line width=0.75pt,black] (0.75,-2.5) to (0.5,-2.75);
    \draw[line width=0.75pt,black] (1.0,-2.5) to (0.75,-2.75);
    \draw[line width=0.75pt,black] (1.25,-2.5) to (1.0,-2.75);
    \draw[line width=0.75pt,black] (1.5,-2.5) to (1.25,-2.75);
    \draw[line width=0.75pt,black] (1.75,-2.5) to (1.5,-2.75);
    \draw[line width=0.75pt,black] (2.0,-2.5) to (1.75,-2.75);
    \draw[line width=0.75pt,black] (2.25,-2.5) to (2.0,-2.75);
    \draw[line width=0.75pt,black] (2.5,-2.5) to (2.25,-2.75);
    \end{tikzpicture}
\caption{SENT experiment}
\label{sec4:fig:tension}
\end{minipage}
\begin{minipage}[b]{0.33\linewidth}
\centering
\begin{tikzpicture}[scale=0.7]
    \coordinate (K) at (0,0);
    % Square
    \draw[line width=0.75pt,black] (-2.5,-2.5) to (2.5,-2.5);
    \draw[line width=0.75pt,black] (2.5,-2.5) to (2.5,2.5);
    \draw[line width=0.75pt,black] (2.5,2.5) to (-2.5,2.5);
    \draw[line width=0.75pt,black] (-2.5,2.5) to (-2.5,-2.5);
    \draw[line width=1.5pt,red] (-2.5,0) to (0,0);
    % Top Boundary
    \draw[line width=0.75pt,black] (2.5,2.75) to (-2.5,2.75);
    \draw[->,line width=1.5pt,black] (0.1,3.25) to (0.8,3.25);
    \node[ ] at (-0.5,3.25) {$\Delta\disp$};
    %\node[ ] at (0,3.85) {$\disp_y=0$};
    % Bottom Boundary
    \draw[line width=0.75pt,black] (-2.5,-2.5) to (-2.75,-2.75);
    \draw[line width=0.75pt,black] (-2.25,-2.5) to (-2.5,-2.75);
    \draw[line width=0.75pt,black] (-2.0,-2.5) to (-2.25,-2.75);
    \draw[line width=0.75pt,black] (-1.75,-2.5) to (-2.0,-2.75);
    \draw[line width=0.75pt,black] (-1.5,-2.5) to (-1.75,-2.75);
    \draw[line width=0.75pt,black] (-1.25,-2.5) to (-1.5,-2.75);
    \draw[line width=0.75pt,black] (-1.0,-2.5) to (-1.25,-2.75);
    \draw[line width=0.75pt,black] (-0.75,-2.5) to (-1.0,-2.75);
    \draw[line width=0.75pt,black] (-0.5,-2.5) to (-0.75,-2.75);
    \draw[line width=0.75pt,black] (-0.25,-2.5) to (-0.5,-2.75);
    \draw[line width=0.75pt,black] (0.0,-2.5) to (-0.25,-2.75);
    \draw[line width=0.75pt,black] (0.25,-2.5) to (0.0,-2.75);
    \draw[line width=0.75pt,black] (0.5,-2.5) to (0.25,-2.75);
    \draw[line width=0.75pt,black] (0.75,-2.5) to (0.5,-2.75);
    \draw[line width=0.75pt,black] (1.0,-2.5) to (0.75,-2.75);
    \draw[line width=0.75pt,black] (1.25,-2.5) to (1.0,-2.75);
    \draw[line width=0.75pt,black] (1.5,-2.5) to (1.25,-2.75);
    \draw[line width=0.75pt,black] (1.75,-2.5) to (1.5,-2.75);
    \draw[line width=0.75pt,black] (2.0,-2.5) to (1.75,-2.75);
    \draw[line width=0.75pt,black] (2.25,-2.5) to (2.0,-2.75);
    \draw[line width=0.75pt,black] (2.5,-2.5) to (2.25,-2.75);
    % Left edge 1
    \draw[fill=gray!50] (-2.75,1.25) -- (-2.25,1.25) -- (-2.5,1.75)-- (-2.75,1.25);
    \draw[fill=black!75] (-2.5,1.05) circle (0.2);
    \draw[line width=1pt,black] (-2.25,0.8) to (-2.75,0.8);
    % Left edge 2
    \draw[fill=gray!50] (-2.75,-1.5) -- (-2.25,-1.5) -- (-2.5,-1.)-- (-2.75,-1.5);
    \draw[fill=black!75] (-2.5,-1.7) circle (0.2);
    \draw[line width=1pt,black] (-2.25,-1.95) to (-2.75,-1.95);
    % Right edge 1
    \draw[fill=gray!50] (2.75,1.25) -- (2.25,1.25) -- (2.5,1.75)-- (2.75,1.25);
    \draw[fill=black!75] (2.5,1.05) circle (0.2);
    \draw[line width=1pt,black] (2.25,0.8) to (2.75,0.8);
    % Right edge 2
    \draw[fill=gray!50] (2.75,-1.5) -- (2.25,-1.5) -- (2.5,-1.)-- (2.75,-1.5);
    \draw[fill=black!75] (2.5,-1.7) circle (0.2);
    \draw[line width=1pt,black] (2.25,-1.95) to (2.75,-1.95);
    \end{tikzpicture}
\caption{SENS experiment}
\label{sec4:fig:miehe_shear}
\end{minipage}
\begin{minipage}[b]{0.33\linewidth}
\centering
\begin{tabular}{ll} \hline
  \textbf{Parameters} & \textbf{Value} \\ \hline
  Fracture Model & AT2 \\
  $E_0$ & 210.0 [GPa] \\
  $\nu$ & 0.3 [-] \\
  $\gc$ & 2.7 [N/mm] \\
  $\l$ & 1.5e-2 [mm] \\
  $\alpha$ & $\beta \gc/\l$ \\ \hline
  \end{tabular}
\captionof{table}{Model parameters}
\label{sec4:table:miehe_tension_shear}
\end{minipage}
\end{figure}

\begin{figure}[!ht]
  \begin{subfigure}[t]{0.45\textwidth}
  \centering
  \begin{tikzpicture}[thick,scale=0.95, every node/.style={scale=1.0}]
    \begin{axis}[width=9.5cm,height=5.75cm,legend style={draw=none,at={(0.5,1.375)},anchor=north,legend cell align=left}, legend columns = 2,
      transpose legend, ylabel={Load\:[kN]},xlabel={Displacement\:[mm]}, xmin=0, ymin=0, xmax=0.0065, ymax=0.8, yticklabel style={/pgf/number format/.cd,fixed,precision=2},
                 every axis plot/.append style={thick}]
    \pgfplotstableread{./Data/SENT/NotchTension_Miehe.txt}\Adata;
    \pgfplotstableread{./Data/SENT/martinez.txt}\Bdata;
    \pgfplotstableread{./Data/SENT/NotchTension_Ambati.txt}\Cdata;
    \pgfplotstableread{./Data/SENT/problem_lodi_10.dat}\Ddata;
    \pgfplotstableread{./Data/SENT/problem_lodi_100.dat}\Edata;
    \pgfplotstableread{./Data/SENT/problem_lodi_250.dat}\Fdata;
    \addplot [ 
           color=black, 
%           only marks, 
           mark=*, 
           mark size=0.05pt, 
         ]
         table
         [
           x expr=\thisrowno{0}, 
           y expr=\thisrowno{1}
         ] {\Adata};
         \addlegendentry{Miehe [39]}
    \addplot [ 
           color=red, 
%           only marks, 
           mark=*, 
           mark size=0.05pt, 
         ]
         table
         [
           x expr=\thisrowno{0}, 
           y expr=\thisrowno{1}/1e3
         ] {\Bdata};
         \addlegendentry{Kristensen [20]}
    \addplot [ 
           color=blue, 
%           only marks, 
           mark=*, 
           mark size=0.25pt, 
         ]
         table
         [
           x expr=\thisrowno{0}, 
           y expr=\thisrowno{1}
         ] {\Cdata};
         \addlegendentry{Ambati [58]} 
    \addplot [ 
           color=green, 
%           only marks, 
           mark=*, 
           mark size=0.05pt, 
         ]
         table
         [
           x expr=\thisrowno{0}, 
           y expr=\thisrowno{1}/1e3
         ] {\Ddata};
         \addlegendentry{10} 
    \addplot [ 
           color=cyan, 
%           only marks, 
           mark=*, 
           mark size=0.05pt, 
         ]
         table
         [
           x expr=\thisrowno{0}, 
           y expr=\thisrowno{1}/1e3
         ] {\Edata};
         \addlegendentry{100}
    \addplot [ 
           color=magenta, 
%           only marks, 
           mark=*, 
           mark size=0.05pt, 
         ]
         table
         [
           x expr=\thisrowno{0}, 
           y expr=\thisrowno{1}/1e3
         ] {\Fdata};
         \addlegendentry{250}
    \end{axis}
    \end{tikzpicture}
    \caption{ }
    \label{sec4:fig:miehe_tension_lodi}
  \end{subfigure}
  \hfill
  \begin{subfigure}[t]{0.45\textwidth}
  \centering
    \begin{tikzpicture}
    \node[inner sep=0pt] () at (0,0)
    {\includegraphics[width=6.0cm,trim=7cm 1cm 7cm 1cm, clip]{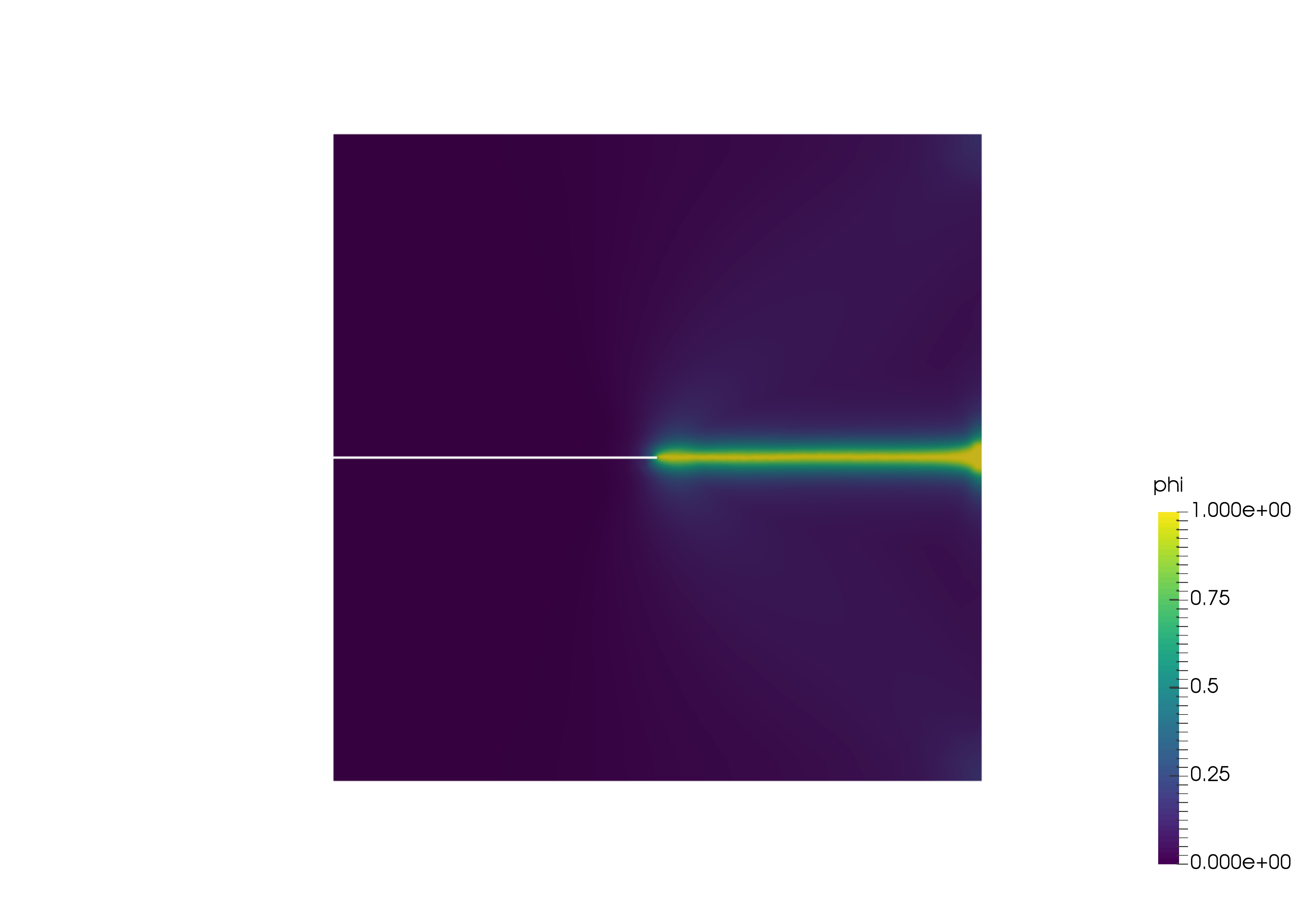}};
    \node[inner sep=0pt] () at (-1.15,-3.05)
    {\begin{axis}[
    hide axis,
    scale only axis,
    height=0pt,
    width=0pt,
    colormap/viridis,
    colorbar horizontal,
    point meta min=0,
    point meta max=1,
    colorbar style={
        width=4.75cm,
        xtick={0,0.5,1.0},
        xticklabel style = {yshift=-0.075cm}
    }]
    \addplot [draw=none] coordinates {(0,0)};
    \end{axis}};
    \node[inner sep=0pt] () at (0,-3.75) {$\pf$};
    \end{tikzpicture}
    \caption{ }
    \label{sec4:fig:pf_sen_tension_failure}
  \end{subfigure}
  \caption{Figure (a) presents the load-displacement curves for the single edge notched specimen under tension. Here, $\beta$ is varied as \{10,100,250\}. Figure (b) shows the distribution of the phase-field variable at the final step of the analysis.}
\end{figure}

\begin{figure}[!ht]
  \begin{subfigure}[t]{0.3\textwidth}
  \centering
  \begin{tikzpicture}[thick,scale=0.95, every node/.style={scale=1.0}]
    \begin{axis}[width=4.75cm,height=5.75cm,legend style={draw=none,at={(0.5,1.375)},anchor=north,legend cell align=left}, legend columns = 2,
      transpose legend, ylabel={Value\:[-]},xlabel={Height\:[mm]}, xmin=0, ymin=0, xmax=1., ymax=1.0, yticklabel style={/pgf/number format/.cd,fixed,precision=2},
                 every axis plot/.append style={thick}]
    \pgfplotstableread[col sep = comma]{./Data/SENT/profile_10.txt}\Adata;
    \addplot [ 
           color=black, 
%           only marks, 
           mark=*, 
           mark size=0.05pt, 
         ]
         table
         [
           x expr=\thisrowno{6}, 
           y expr=\thisrowno{4}
         ] {\Adata};
         \addlegendentry{$\pf$}
    \addplot [ 
           color=red, 
%           only marks, 
           mark=*, 
           mark size=0.05pt, 
         ]
         table
         [
           x expr=\thisrowno{6}, 
           y expr=\thisrowno{3}
         ] {\Adata};
         \addlegendentry{\mpf}
    \end{axis}
    \end{tikzpicture}
    \caption{$\beta = 10$}
    \label{sec4:fig:sentProfile10}
  \end{subfigure}
  \hfill
  \begin{subfigure}[t]{0.3\textwidth}
  \centering
  \begin{tikzpicture}[thick,scale=0.95, every node/.style={scale=1.0}]
    \begin{axis}[width=4.75cm,height=5.75cm,legend style={draw=none,at={(0.5,1.375)},anchor=north,legend cell align=left}, legend columns = 2,
      transpose legend, xlabel={Height\:[mm]}, xmin=0, ymin=0, xmax=1., ymax=1.0, yticklabel style={/pgf/number format/.cd,fixed,precision=2},
                 every axis plot/.append style={thick}]
    \pgfplotstableread[col sep = comma]{./Data/SENT/profile_100.txt}\Adata;
    \addplot [ 
           color=black, 
%           only marks, 
           mark=*, 
           mark size=0.05pt, 
         ]
         table
         [
           x expr=\thisrowno{6}, 
           y expr=\thisrowno{4}
         ] {\Adata};
         \addlegendentry{$\pf$}
    \addplot [ 
           color=red, 
%           only marks, 
           mark=*, 
           mark size=0.05pt, 
         ]
         table
         [
           x expr=\thisrowno{6}, 
           y expr=\thisrowno{3}
         ] {\Adata};
         \addlegendentry{\mpf}
    \end{axis}
    \end{tikzpicture}
    \caption{$\beta = 100$}
    \label{sec4:fig:sentProfile100}
  \end{subfigure}
  \hfill
  \begin{subfigure}[t]{0.3\textwidth}
  \centering
  \begin{tikzpicture}[thick,scale=0.95, every node/.style={scale=1.0}]
    \begin{axis}[width=4.75cm,height=5.75cm,legend style={draw=none,at={(0.5,1.375)},anchor=north,legend cell align=left}, legend columns = 2,
      transpose legend, xlabel={Height\:[mm]}, xmin=0, ymin=0, xmax=1., ymax=1.0, yticklabel style={/pgf/number format/.cd,fixed,precision=2},
                 every axis plot/.append style={thick}]
    \pgfplotstableread[col sep = comma]{./Data/SENT/profile_250.txt}\Adata;
    \addplot [ 
           color=black, 
%           only marks, 
           mark=*, 
           mark size=0.05pt, 
         ]
         table
         [
           x expr=\thisrowno{6}, 
           y expr=\thisrowno{4}
         ] {\Adata};
         \addlegendentry{$\pf$}
    \addplot [ 
           color=red, 
%           only marks, 
           mark=*, 
           mark size=0.05pt, 
         ]
         table
         [
           x expr=\thisrowno{6}, 
           y expr=\thisrowno{3}
         ] {\Adata};
         \addlegendentry{\mpf}
    \end{axis}
    \end{tikzpicture}
    \caption{$\beta = 250$}
    \label{sec4:fig:sentProfile250}
  \end{subfigure}
  \caption{Figures present the phase-field ($\pf$) and the micromorphic variable ($\mpf$) for different $\beta$ values at $x=0.75$ [mm], along the height of the SENT specimen.}
\end{figure}
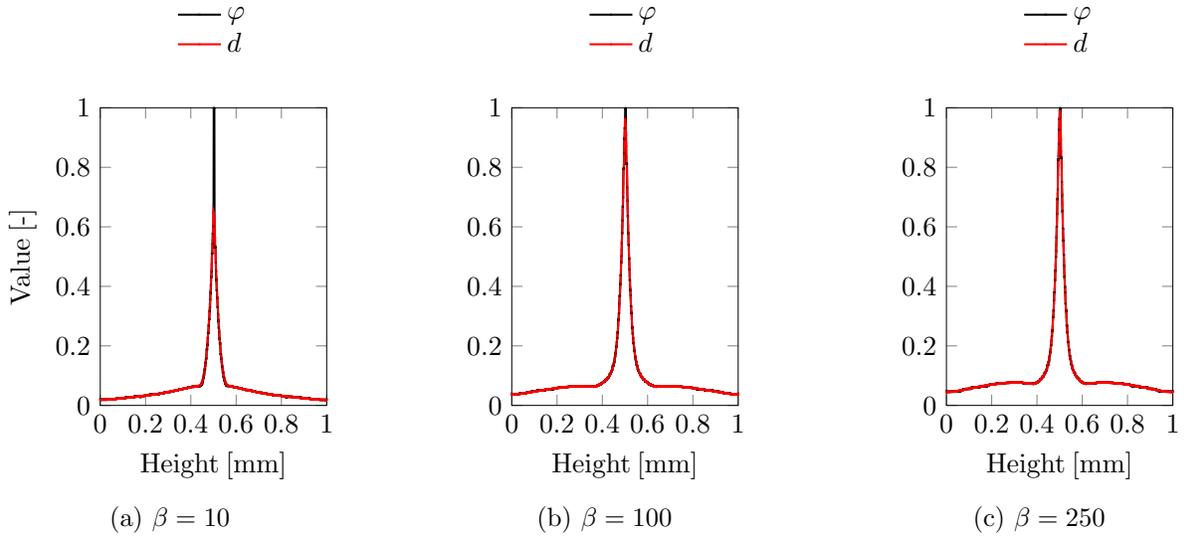

Figure \ref{sec4:fig:miehe_tension_lodi} present the load-displacement curves, obtained using the micromorphic phase-field fracture model with different values of $\beta$. They are compared with the load-displacement curves from the literature \cite{Miehe2010,Ambati2015,KRISTENSEN2020104093}. While $\beta = 100, 250$ yield curves within the range of those from the literature, $\beta = 10$ predicts $\approx 24 \%$ lower peak load. This is accompanied by a lower displacement at failure of the specimen, evident from the green curve in Figure \ref{sec4:fig:miehe_tension_lodi}. In order to provide an explanation for this behaviour, both the phase-field and the micrmorphic variable are plotted across the section, $x = 0.75$ [mm] in Figures \ref{sec4:fig:sentProfile10}-\ref{sec4:fig:sentProfile250}, for varying $\beta$. For $\beta = 10$, the phase-field ($\pf$) and the micromorphic variable ($\mpf$) are not identical, indicating an insufficient interaction between the local and the non-local fields. Besides, comparing Figure \ref{sec4:fig:sentProfile10} with Figures \ref{sec4:fig:sentProfile100} and \ref{sec4:fig:sentProfile250}, a second observation maybe made w.r.t. to the regularization. For $\beta = 10$, the phase-field values higher than $0.7$ are not regularized. This results in local material response as the strength parameters (stress, stiffness) are depedent on the phase-field. This explains the different load-displacement response compared to those obtained from $\beta = 100, 250$.
For the latter values, the phase-field fracture topology in the final step of the simulation in Figure \ref{sec4:fig:pf_sen_tension_failure}, is similar to those from the literature \cite{Miehe2010,Ambati2015,KRISTENSEN2020104093}.

% Profiles

\subsection{Single Edge Notched specimen under Shear (SENS)}\label{sec4:SENshear}

In order to perform a shear test, the single edge notched specimen is loaded horizontally along the top edge as shown in Figure \ref{sec4:fig:miehe_shear}. The material properties remain same as presented in Table \ref{sec4:table:miehe_tension_shear}. A quasi-static loading is applied to the top boundary in the form of prescribed displacement increment $\Delta\disp = 1e-4$[mm] for the first $85$ steps, following which it is changed to $5e-6$[mm]. Furthermore, the bottom boundary remains fixed, and roller support is implemented in left and right edges thereby restricting the vertical displacement.

\begin{figure}[!ht]
  \begin{subfigure}[t]{0.45\textwidth}
  \centering
  \begin{tikzpicture}[thick,scale=0.95, every node/.style={scale=1.0}]
    \begin{axis}[width=9.5cm,height=5.75cm,legend style={draw=none,at={(0.5,1.375)},anchor=north,legend cell align=left}, legend columns = 3,
       ylabel={Load\:[kN]},xlabel={Displacement\:[mm]}, xmin=0, ymin=0, xmax=0.013, ymax=0.55, yticklabel style={/pgf/number format/.cd,fixed,precision=2},
                 every axis plot/.append style={thick}]
    \pgfplotstableread{./Data/SENS/problem_lodi100.dat}\Adata;
    \pgfplotstableread{./Data/SENS/problem_lodi250.dat}\Cdata;
    \pgfplotstableread{./Data/SENS/problem_lodi10.dat}\Ddata;
    \pgfplotstableread[col sep = comma]{./Data/SENS/NotchShear_Wu.txt}\Bdata;
    \addplot [ 
           color=black, 
%           only marks, 
           mark=*, 
           mark size=0.05pt, 
         ]
         table
         [
           x expr=\thisrowno{0}, 
           y expr=\thisrowno{1}
         ] {\Bdata};
         \addlegendentry{Wu et. al. [59]}
    \addplot [ 
           color=green, 
%           only marks, 
           mark=*, 
           mark size=0.05pt, 
         ]
         table
         [
           x expr=\thisrowno{0}, 
           y expr=\thisrowno{1}/1e3
         ] {\Ddata};
         \addlegendentry{10}
    \addplot [ 
           color=cyan, 
%           only marks, 
           mark=*, 
           mark size=0.05pt, 
         ]
         table
         [
           x expr=\thisrowno{0}, 
           y expr=\thisrowno{1}/1e3
         ] {\Adata};
         \addlegendentry{100}
    \addplot [ 
           color=magenta, 
%           only marks, 
           mark=*, 
           mark size=0.05pt, 
         ]
         table
         [
           x expr=\thisrowno{0}, 
           y expr=\thisrowno{1}/1e3
         ] {\Cdata};
         \addlegendentry{250}
    \end{axis}
    \end{tikzpicture}
    \caption{ }
    \label{sec4:fig:miehe_shear_lodi}
  \end{subfigure}
  \hfill
  \begin{subfigure}[t]{0.45\textwidth}
  \centering
    \begin{tikzpicture}
    \node[inner sep=0pt] () at (0,0)
    {\includegraphics[width=6.0cm,trim=7cm 1cm 7cm 1cm, clip]{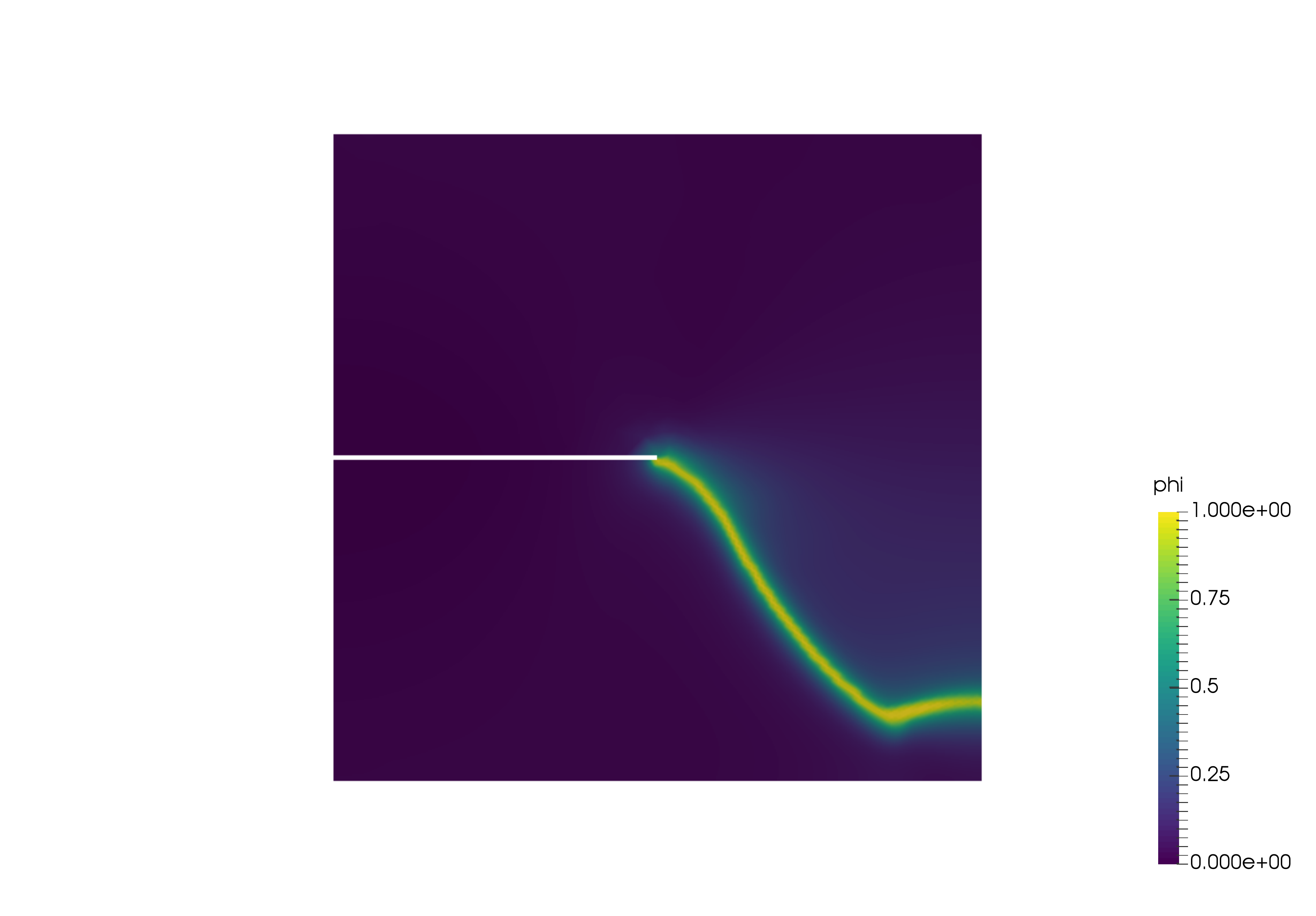}};
    \node[inner sep=0pt] () at (-1.15,-3.05)
    {\begin{axis}[
    hide axis,
    scale only axis,
    height=0pt,
    width=0pt,
    colormap/viridis,
    colorbar horizontal,
    point meta min=0,
    point meta max=1,
    colorbar style={
        width=4.75cm,
        xtick={0,0.5,1.0},
        xticklabel style = {yshift=-0.075cm}
    }]
    \addplot [draw=none] coordinates {(0,0)};
    \end{axis}};
    \node[inner sep=0pt] () at (0,-3.75) {$\pf$};
    \end{tikzpicture}
    \caption{ }
    \label{sec4:fig:pf_sen_shear_failure}
  \end{subfigure}
  \caption{Figure (a) presents the load-displacement curves for the single edge notched specimen under shear. Here, $\beta$ is varied as \{100,250\}. Figure (b) shows the distribution of the phase-field variable at the final step of the analysis.}
\end{figure}

Figure \ref{sec4:fig:miehe_shear_lodi} shows the load displacement curves obtained using $\beta = 10, 100, 250$. The curves obtained using $\beta = 100, 250$ is similar to that obtained in \cite{Wu20201review}. Moreoever, the phase-field topology at the final step of the analysis (see Figure \ref{sec4:fig:pf_sen_shear_failure}) is similar to that presented in Figure 5 of \cite{Wu20201review}. For $\beta = 10$, the non-local behaviour of the micromorphic model is not sufficient. This exaplains the loss of integrity of specimen at an earlier stage in the loading process (cf. Figures \ref{sec4:fig:miehe_shear_lodi} and \ref{sec4:fig:miehe_tension_lodi}).

\subsection{Winkler L-panel}\label{sec4:Lpanel}

The concrete L-shaped panel studied by \cite{winkler2001traglastuntersuchungen,unger2007modelling} is considered in this sub-section. Figure \ref{sec4:fig:Lpanel} shows the geometry as well as the loading conditions. The longer edges of the panel are 500 [mm] and the smaller edges are 250 [mm]. The loading is applied on the edge marked in blue, 30 [mm] in length, and is in the form of displacement increments of $\Delta\disp = 1e-3$ [mm]. The additional model parameters are presented in Table \ref{sec4:table:Lpanel_params}.

\begin{figure}[ht]
\begin{minipage}[b]{0.45\linewidth}
\centering
\begin{tikzpicture}[scale=0.6]
    \coordinate (K) at (0,0);
    % Square
    \draw[line width=0.75pt,black] (-4,-4) to (0,-4);
    \draw[line width=0.75pt,black] (0,-4) to (0,0);
    \draw[line width=0.75pt,black] (0,0) to (4,0);
    \draw[line width=0.75pt,black] (4,0)  to (4,4);
    \draw[line width=0.75pt,black] (4,4) to (-4,4);
    \draw[line width=0.75pt,black] (-4,4) to (-4,-4);
    % Loading
    \draw[line width=2.0pt,blue] (3.5,0) to (4,0);
    \draw[->,line width=1.5pt,black] (3.75,-0.9) to (3.75,-0.25);
    \node[ ] at (3.75,-1.15) {$\Delta\disp$};
    % Left edge
    %\draw[fill=black!75] (-2.5,1.05) circle (0.2);
    \draw[line width=1pt,black] (0,-4) to (-0.25,-4.25);
    \draw[line width=1pt,black] (-0.25,-4.) to (-0.5,-4.25);
    \draw[line width=1pt,black] (-0.5,-4.) to (-0.75,-4.25);
    \draw[line width=1pt,black] (-0.75,-4.) to (-1,-4.25);
    \draw[line width=1pt,black] (-1,-4.) to (-1.25,-4.25);
    \draw[line width=1pt,black] (-1.25,-4.) to (-1.5,-4.25);
    \draw[line width=1pt,black] (-1.5,-4.) to (-1.75,-4.25);
    \draw[line width=1pt,black] (-1.75,-4.) to (-2.,-4.25);
    \draw[line width=1pt,black] (-2.,-4.) to (-2.25,-4.25);
    \draw[line width=1pt,black] (-2.25,-4.) to (-2.5,-4.25);
    \draw[line width=1pt,black] (-2.5,-4.) to (-2.75,-4.25);
    \draw[line width=1pt,black] (-2.75,-4.) to (-3.,-4.25);
    \draw[line width=1pt,black] (-3.,-4.) to (-3.25,-4.25);
    \draw[line width=1pt,black] (-3.25,-4.) to (-3.5,-4.25);
    \draw[line width=1pt,black] (-3.5,-4.) to (-3.75,-4.25);
    \draw[line width=1pt,black] (-3.75,-4.) to (-4.,-4.25);
    \draw[line width=1pt,black] (-4.,-4.) to (-4.25,-4.25);
    \end{tikzpicture}
\caption{Winkler L-panel}
\label{sec4:fig:Lpanel}
\end{minipage}
\hspace{0.5cm}
\begin{minipage}[b]{0.45\linewidth}
\centering
\begin{tabular}{ll} \hline
  \textbf{Parameters}  & \textbf{Value} \\ \hline
  Model & Quasi-brittle \\
  Softening & Cornellisen et. al. \cite{cornelissen1986experimental} \\
  $E_0$  & 2.0e4 [MPa] \\
  $\nu$  & 0.18 [-] \\
  $f_t$  & 2.5 [MPa] \\
  $\gc$  & 0.130 [N/mm] \\
  $\l$  & 10 [mm] \\
  $\alpha$ & $\beta \gc/\l$ \\ \hline
  \end{tabular}
\captionof{table}{Parameters for L-shaped panel test \cite{unger2007modelling}}
\label{sec4:table:Lpanel_params}
\end{minipage}
\end{figure}

\begin{figure}[!ht]
\centering
  \begin{subfigure}[t]{0.45\textwidth}
  \centering
    \begin{tikzpicture}[thick,scale=0.95, every node/.style={scale=1.175}]
    \begin{axis}[legend style={draw=none}, legend columns = 3,
      transpose legend, ylabel={Load\:[kN]},xlabel={Displacement\:[mm]}, xmin=0, ymin=0, xmax=1.0, ymax=8.0, yticklabel style={
        /pgf/number format/fixed,
        /pgf/number format/precision=5
        },
        scaled y ticks=false,
        xticklabel style={
        /pgf/number format/fixed,
        /pgf/number format/precision=5
        },
        scaled x ticks=false,
        every axis plot/.append style={thick}]
    \pgfplotstableread{./Data/WinklerL/problem_lodi250.dat}\Edata;
    \pgfplotstableread{./Data/WinklerL/problem_lodi100.dat}\Cdata;
    \pgfplotstableread[col sep = comma]{./Data/WinklerL/LPanelUpperLodi.txt}\Bdata;
    \pgfplotstableread[col sep = comma]{./Data/WinklerL/LPanelLowerLodi.txt}\Adata;
    \addplot [ name path = A,
           color=blue!20, 
%           only marks, 
           mark=*, 
           mark size=0.01pt,
           forget plot
         ]
         table
         [
           x expr=\thisrowno{0}, 
           y expr=\thisrowno{1}
         ] {\Adata};
    \addplot [ name path = B,
           color=blue!20, 
%           only marks, 
           mark=*, 
           mark size=0.01pt, 
           forget plot
         ]
         table
         [
           x expr=\thisrowno{0}, 
           y expr=\thisrowno{1}
         ] {\Bdata};
    % Fill area between paths
    \addplot [blue!20, forget plot] fill between [of = A and B];     
    \addplot [ 
           color=cyan, 
%           only marks, 
           mark=*, 
           mark size=0.05pt, 
         ]
         table
         [
           x expr=\thisrowno{0}, 
           y expr=\thisrowno{1}/10
         ] {\Cdata};
         \addlegendentry{100}
    \addplot [ 
           color=magenta, 
%           only marks, 
           mark=*, 
           mark size=0.05pt, 
         ]
         table
         [
           x expr=\thisrowno{0}, 
           y expr=\thisrowno{1}/10
         ] {\Edata};
         \addlegendentry{250}     
    \end{axis}
    \end{tikzpicture}
    \caption{Load-displacement plot}
    \label{sec4:fig:LPanel_lodi}
  \end{subfigure}
  \hspace{5mm}
  \begin{subfigure}[t]{0.45\textwidth}
  \centering
    \begin{tikzpicture}
    \node[inner sep=0pt] () at (0,1.5)
    {\includegraphics[width=6.0cm,trim=7cm 1cm 7cm 1cm, clip]{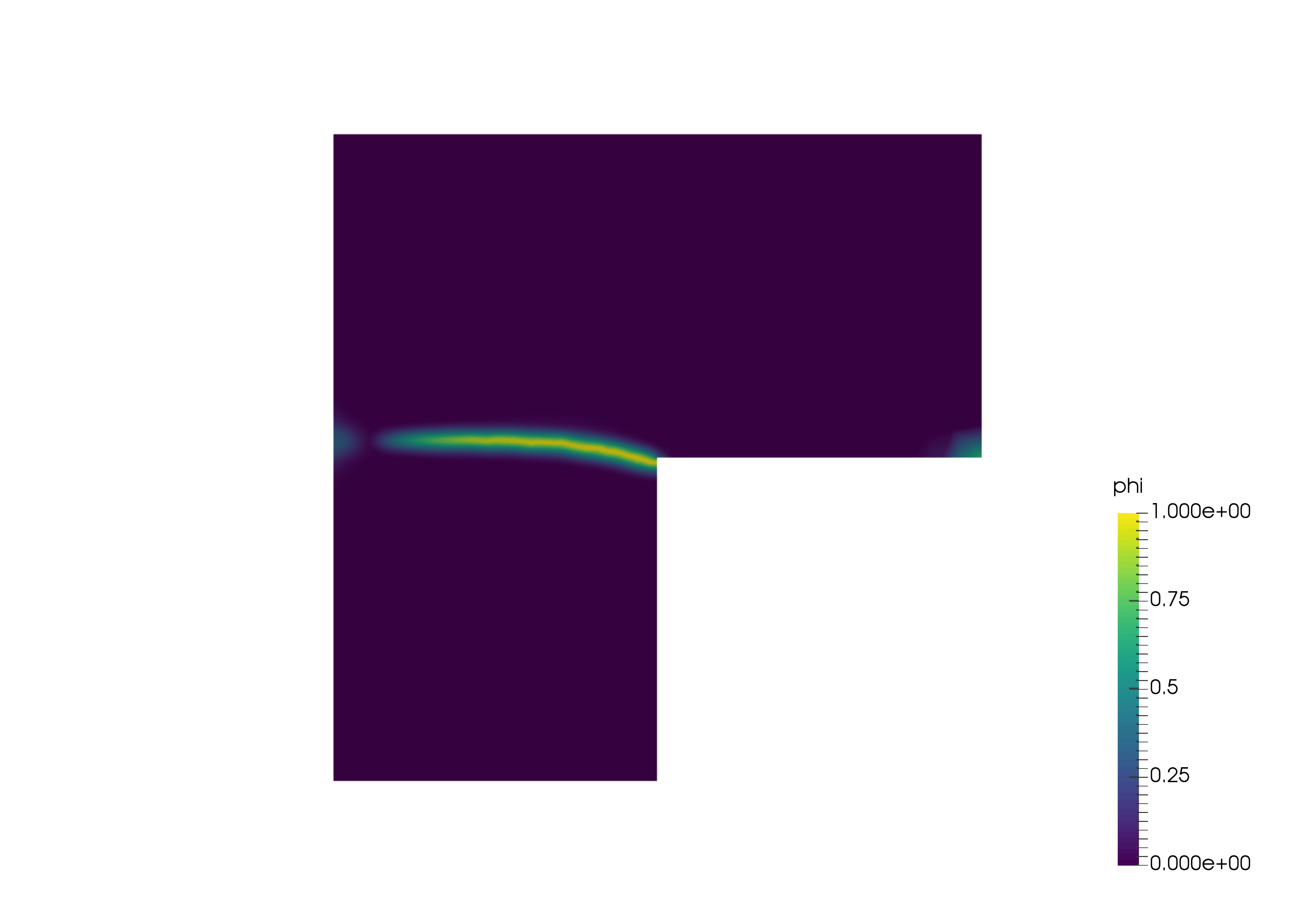}};
    \node[inner sep=0pt] () at (-1.2,-1.75)
    {\begin{axis}[
    hide axis,
    scale only axis,
    height=0pt,
    width=0pt,
    colormap/viridis,
    colorbar horizontal,
    point meta min=0,
    point meta max=1,
    colorbar style={
        width=4.75cm,
        xtick={0,0.5,1.0},
        xticklabel style = {yshift=-0.075cm}
    }]
    \addplot [draw=none] coordinates {(0,0)};
    \end{axis}};
    \node[inner sep=0pt] () at (0,-2.45) {$\pf$};
    \end{tikzpicture}
    \caption{ }
    \label{sec4:fig:pf_LPanel_failure}
  \end{subfigure}
  \caption{Figure (a) presents the load-displacement curves for the concrete Winkler L-panel test, for $\beta = 100, 250$. The experimental range is represented by the shaded area. Figure (b) shows the distribution of the phase-field variable at the final step of the analysis.}
  \label{sec4:fig:LPanel_results}
\end{figure}

Figure \ref{sec4:fig:LPanel_lodi} presents the load-displacement curves using the micromorphic phase-field fracture model, with $\beta = 100,250$. The range of the experimentally obtained load-displacement curves is shown in the shaded region. For both values of $\beta$, the load-displacement curves exhibit a good agreement with the experimental region. Moreover, the phase-field topology in the final step of the simulation in Figure \ref{sec4:fig:pf_LPanel_failure} is similar to that from the literature \cite{winkler2001traglastuntersuchungen}.

\subsection{Quasi-brittle: Concrete three-point bending}\label{sec4:concreteTPB}

A three-point bending experiment on a notched concrete beam reported in \cite{Rots1988} is considered here. The beam has dimensions $450\times100$ [mm$^2$], and has a notch $5\times50$ [mm$^2$]. A schematic of the beam along with the loading conditions is presented in Figure \ref{sec4:fig:3pt_bending}. Displacement-based load increments of $\Delta\disp = 1e-3$ [mm] is enforced throughout the simulation. The model parameters are presented in Table \ref{sec4:table:3pt_bending}.

\begin{figure}[ht]
\begin{minipage}[b]{0.45\linewidth}
\centering
\begin{tikzpicture}[scale=0.7]
    \coordinate (K) at (0,0);
    % Square
    \draw[line width=0.75pt,black] (4,2) to (-4,2);
    \draw[line width=0.75pt,black] (4,2) to (4,0);
    \draw[line width=0.75pt,black] (4,0) to (0.05,0);
    \draw[line width=0.75pt,black] (0.05,0) to (0.05,1.0);
    \draw[line width=0.75pt,black] (0.05,1.0) to (-0.05,1.0);
    \draw[line width=0.75pt,black] (-0.05,1.0) to (-0.05,0.0);
    \draw[line width=0.75pt,black] (-0.05,0.0) to (-4.0,0.0);
    \draw[line width=0.75pt,black] (-4.0,0.0) to (-4.0,2);
    % Loading
    \draw[->,line width=1.5pt,black] (0.0,2.5) to (0.0,2.);
    \node[ ] at (0,2.75) {$\Delta\disp$};
    % Left edge
    \draw[fill=gray!50] (-4.25,-0.5) -- (-3.75,-0.5) -- (-4,0)-- (-4.25,-0.5);
    %\draw[fill=black!75] (-2.5,1.05) circle (0.2);
    \draw[line width=1pt,black] (-4.5,-0.5) to (-3.5,-0.5);
    \draw[line width=1pt,black] (-4.5,-0.5) to (-4.75,-0.75);
    \draw[line width=1pt,black] (-4.25,-0.5) to (-4.5,-0.75);
    \draw[line width=1pt,black] (-4,-0.5) to (-4.25,-0.75);
    \draw[line width=1pt,black] (-3.75,-0.5) to (-4,-0.75);
    \draw[line width=1pt,black] (-3.5,-0.5) to (-3.75,-0.75);
    % Right edge
    \draw[fill=gray!50] (4.25,-0.5) -- (3.75,-0.5) -- (4,0)-- (4.25,-0.5);
    \draw[fill=black!75] (4,-0.7) circle (0.2);
    \draw[line width=1pt,black] (4.5,-0.95) to (3.5,-0.95);
    \end{tikzpicture}
\caption{Three point bending test}
\label{sec4:fig:3pt_bending}
\end{minipage}
\hspace{0.5cm}
\begin{minipage}[b]{0.45\linewidth}
\centering
\begin{tabular}{ll} \hline
  \textbf{Parameters}  & \textbf{Value} \\ \hline
  Model & Quasi-brittle \\
  Softening & Cornellisen et. al. \cite{cornelissen1986experimental} \\
  $E_0$  & 2e4 [MPa] \\
  $\nu$  & 0.2 [-] \\
  $f_t$  & 2.4 [MPa] \\
  $\gc$  & 0.113 [N/mm] \\
  $\l$  & 2.5 [mm] \\
  $\alpha$ & $\beta \gc/\l$ \\ \hline
  \end{tabular}
\captionof{table}{Parameters for three point bending test}
\label{sec4:table:3pt_bending}
\end{minipage}
\end{figure}

\begin{figure}[!ht]
\centering
  \begin{subfigure}[t]{0.45\textwidth}
  \centering
    \begin{tikzpicture}[thick,scale=0.95, every node/.style={scale=1.175}]
    \begin{axis}[legend style={draw=none}, legend columns = 3,
      transpose legend, ylabel={Load\:[kN]},xlabel={Displacement\:[mm]}, xmin=0, ymin=0, xmax=1.0, ymax=1.75, yticklabel style={
        /pgf/number format/fixed,
        /pgf/number format/precision=5
        },
        scaled y ticks=false,
        xticklabel style={
        /pgf/number format/fixed,
        /pgf/number format/precision=5
        },
        scaled x ticks=false,
        every axis plot/.append style={very thick}]
    \pgfplotstableread{./Data/Rots/problem_lodi_250.dat}\Edata;
    \pgfplotstableread{./Data/Rots/problem_lodi_100.dat}\Cdata;
    \pgfplotstableread[col sep = comma]{./Data/Rots/RotsTPBUpperLodi.txt}\Bdata;
    \pgfplotstableread[col sep = comma]{./Data/Rots/RotsTPBLowerLodi.txt}\Adata;
    \addplot [name path = A, 
           color=blue!20, 
%           only marks, 
           mark=*, 
           mark size=0.25pt, 
           forget plot
         ]
         table
         [
           x expr=\thisrowno{0}, 
           y expr=\thisrowno{1}
         ] {\Adata};
    \addplot [name path = B,
           color=blue!20, 
%           only marks, 
           mark=*, 
           mark size=0.25pt, 
           forget plot
         ]
         table
         [
           x expr=\thisrowno{0}, 
           y expr=\thisrowno{1}
         ] {\Bdata};
    \addplot [blue!20, forget plot] fill between [of = A and B]; 
    \addplot [ 
           color=cyan, 
%           only marks, 
           mark=*, 
           mark size=0.05pt, 
         ]
         table
         [
           x expr=\thisrowno{0}, 
           y expr=\thisrowno{1}
         ] {\Cdata};
         \addlegendentry{100}
    \addplot [ 
           color=magenta, 
%           only marks, 
           mark=*, 
           mark size=0.05pt, 
         ]
         table
         [
           x expr=\thisrowno{0}, 
           y expr=\thisrowno{1}
         ] {\Edata};
         \addlegendentry{250} 
    \end{axis}
    \end{tikzpicture}
    \caption{Load-displacement plot}
    \label{sec4:fig:3pt_bending_lodi}
  \end{subfigure}
  \hspace{5mm}
  \begin{subfigure}[t]{0.45\textwidth}
  \centering
    \begin{tikzpicture}
    \node[inner sep=0pt] () at (0,0)
    {\includegraphics[width=6.0cm,trim=7cm 1cm 7cm 1cm, clip]{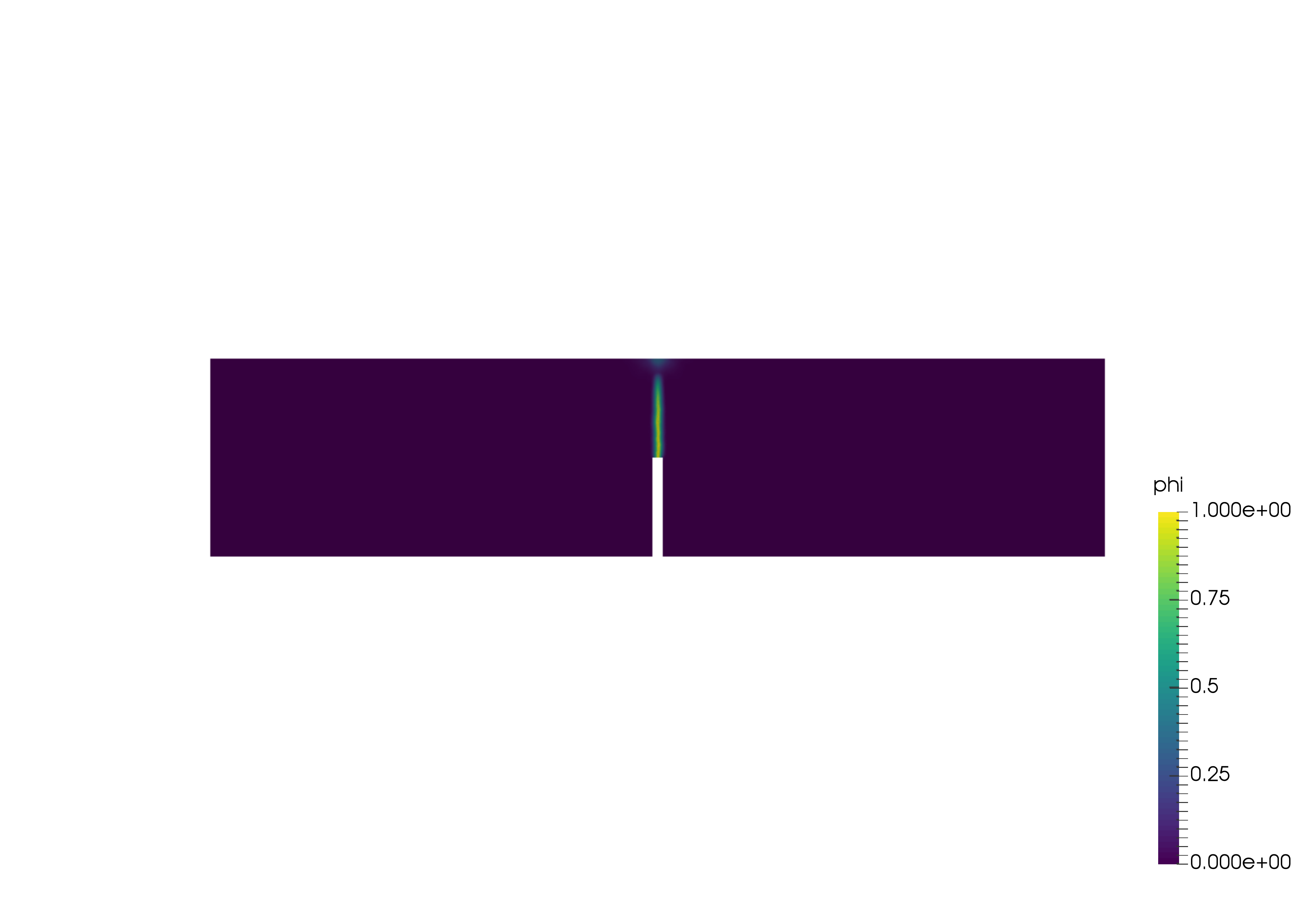}};
    \node[inner sep=0pt] () at (-1.2,-1.75)
    {\begin{axis}[
    hide axis,
    scale only axis,
    height=0pt,
    width=0pt,
    colormap/viridis,
    colorbar horizontal,
    point meta min=0,
    point meta max=1,
    colorbar style={
        width=4.75cm,
        xtick={0,0.5,1.0},
        xticklabel style = {yshift=-0.075cm}
    }]
    \addplot [draw=none] coordinates {(0,0)};
    \end{axis}};
    \node[inner sep=0pt] () at (0,-2.45) {$\pf$};
    \end{tikzpicture}
    \caption{ }
    \label{sec4:fig:pf_3pt_bending_failure}
  \end{subfigure}
  \caption{Figure (a) presents the load-displacement curves for the concrete three-point bending test, with $\beta = 100, 250$. The experimental range is represented by the shaded area. Figure (b) shows the distribution of the phase-field variable at the final step of the analysis in a section of the beam.}
  \label{sec4:fig:3pt_bending_results}
\end{figure}

Figure \ref{sec4:fig:3pt_bending_lodi} presents the load-displacement curves using the micromorphic phase-field fracture model, with $\beta = 100,250$. The range of the experimentally obtained load-displacement curves is exhibited by the shaded region. For both values of $\beta$, the load-displacement curves exhibit a good agreement with the experimental region. Moreover, the phase-field topology in the final step of the simulation in Figure \ref{sec4:fig:pf_3pt_bending_failure} is similar to that from the literature \cite{Rots1988}.

\section{Conclusions and Outlook}\label{sec5}

A novel phase-field fracture model is proposed in this manuscript, based on the micromorhic extension of the phase-field fracture energy functional. In this model, the phase-field variable is local, and a micromorphic variable is introduced for regularization. Combined with an irreversibility criterion on the phase-field variable, this results in a local variational inequality problem. The problem is then trivially solved pointwise (at integration points) throughout the computational domain. For brittle AT1 and AT2 models, an explicit expression for the local phase-field exists, however, for quasi-brittle fracture models, a nonlinear scalar equation needs to be solved. Furthermore, the local nature of the phase-field variable also provides the ease in implementing bounds, i.e., $\pf \in [0,1]$ using the trivial `\textit{min}' and `\textit{max}' operations. Unlike the penalization technique proposed in \cite{GERASIMOV2019990}, fracture irreversibility is achieved with system precision without the need for a high penalty parameter. As regards to computational effort, since the phase-field variable is solved locally, the size of the finite element problem is identical to that of a standard phase-field implementation.

Furthermore, the micromorphic phase-field fracture model does converge towards the original phase-field model in the terms of the energy functional equivalence, when the interaction parameter $\alpha$ is chosen such that $\pf \approx \mpf$. This has been demonstrated in this manuscript with numerical experiments on the single edge notched specimen loaded in tension. When the interaction parameter is sufficiently high ($\approx 250 \gc/\l$), the phase-field and the micromorphic variable profiles are similar for arbitrary sections in a computational domain. This essentially recovers the original phase-field fracture energy functional, while enabling the ease of local phase-field fracture irreversibility. Numerical experiments were also carried out in the context of quasi-brittle fracture, the Winkler L-shaped panel and the concrete three-point bending tests. It is observed that the micromorphic phase-field fracture model results in fracture topologies and load-displacement curves within the experimentally observed range. This demonstrates capability of the micromorphic approach to demonstrate both brittle and quasi-brittle phenomenon. 

Finally, the novel micromorphic phase-field fracture model opens a plethora of future research extensions, particularly, in multi-physics applications, composite laminates \cite{BUI2021107705}. Other studies may include the implementation of a dissipation-based arc-length method \cite{may2015numerical,BHARALI2022114927} or quasi-Newton methods \cite{KRISTENSEN2020102446,WU2020112704} for addressing the non-convexity of energy functional. 

\section{Software Implementation and Data Availability}

The numerical study in Section \ref{sec5} is carried out using an in-house finite element code, inspired by ofeFrac \cite{NGUYENTHANH2020102925}, and based on the Jem and Jive C++ libraries from the Dynaflow Research Group, the Netherlands. The simulation data would be made available in Github repository of the corresponding author (\url{https://github.com/ritukeshbharali}).

\section*{Declaration of Competing Interest}

The authors declare that they have no known competing financial interests or personal relationships that could have appeared to influence the work reported in this paper.

\section*{Acknowledgements}

The financial support from the Swedish Research Council for Sustainable Development (FORMAS) under Grant 2018-01249 and the Swedish Research Council (VR) under Grant 2017-05192 is gratefully acknowledged. The first author would also like to thank Dr. Erik Jan Lingen at Dynaflow Research Group for support with the Jem and Jive libraries.

\begin{appendices}
\numberwithin{equation}{section}
\numberwithin{figure}{section}
\section{Constitutive relation with Amor split \cite{pham2011}}\label{appendixA}

In this manuscript, the fracture driving and residual strain energy densities are computed using the Amor split \cite{pham2011}. With this approach, the fracture is driven by the deviatoric part of the strain energy density together with a positive volumetric part of the strain energy density. The negative volumetric part corresponds to the residual energy. Adopting a Voigt notation, the constitutive material stiffness $\mathbf{D}$ in (\ref{eq:stiff_comp},\ref{eq:stiff_comp2}),

\begin{equation}
    \mathbf{D} = g(\hat{\pf}) \dfrac{\partial \stress^+}{\partial \strain} + \dfrac{\partial \stress^-}{\partial \strain}
\end{equation}

\noindent is computed as,

\begin{equation}
    \mathbf{D} = g(\hat{\pf}) \bigg( K \mathcal{H}(tr(\strain)) \mathbb{P}_{vol} + 2 \mu \mathbb{P}_{dev} \bigg) + K \mathcal{H}(-tr(\strain)) \mathbb{P}_{vol}
\end{equation}

\noindent with bulk and shear modulus $K,\mu$, heaviside function $\mathcal{H}$ and projection matrices,

\begin{equation}
    \mathbb{P}_{vol} = \begin{bmatrix}
1 & 1  & 1  & 0 & 0 & 0 \\ 
1 & 1  & 1  & 0 & 0 & 0 \\ 
1 & 1  & 1  & 0 & 0 & 0 \\ 
0 & 0 & 0 & 0 & 0 & 0\\ 
0 & 0 & 0 & 0 & 0 & 0\\ 
0 & 0 & 0 & 0 & 0 & 0
\end{bmatrix}
\end{equation}

\noindent and

\begin{equation}
    \mathbb{P}_{dev} = \begin{bmatrix}
2/3 & -1/3  & 1/3  & 0 & 0 & 0 \\ 
1/3 & 2/3  & 1/3  & 0 & 0 & 0 \\ 
1/3 & 1/3  & 2/3  & 0 & 0 & 0 \\ 
0 & 0 & 0 & 1/2 & 0 & 0\\ 
0 & 0 & 0 & 0 & 1/2 & 0\\ 
0 & 0 & 0 & 0 & 0 & 1/2
\end{bmatrix}.
\end{equation}

\end{appendices}

\printbibliography

\end{document}